\documentclass{cocv}
\usepackage{graphicx}
\usepackage{wrapfig}
\begin{document}
\title{Method for finding solution to \\ ``quasidifferentiable'' differential inclusion}
\thanks{The main results of this paper (Sections 4, 5, 6, 7) were obtained in Saint Petersburg State University and supported by Russian Science Foundation (grant  23-71-01065).}
\author{Alexander Fominyh}\address{St.\,Petersburg University, 7--9, Universitetskaya nab., St.\,Petersburg, 199034, Russian Federation}
%
%
\begin{abstract} The paper explores the differential inclusion of a special form. It is supposed that the support function of the set in the right-hand side of an inclusion may contain the sum of the maximum and the minimum of the finite number of continuously differentiable (in phase coordinates) functions. It is required to find a trajectory that would satisfy differential inclusion with the boundary conditions prescribed and simultaneously lie on the surface given. We give substantial examples of problems where such differential inclusions may occur: models of discontinuous systems, linear control systems where the control function or/and disturbance of the right-hand side is/are known to be subject to some nonsmooth (in phase vector) constraints, some real mechanical models and differential inclusions per se with special geometrical structure of the right-hand side. The initial problem is reduced to a variational one. It is proved that the resulting functional to be minimized is quasidifferentiable. The necessary minimum conditions in terms of quasidifferential are formulated. The steepest (or the quasidifferential) descent method in a classical form is then applied to find stationary points of the functional obtained. Herewith, the functional is constructed in such a way that one can verify whether the stationary point constructed is indeed a global minimum point of the problem. The ``weak'' convergence of the method proposed is proved for some particular cases. The method constructed is illustrated by numerical examples.  \end{abstract}
%
%
\maketitle
\section{Introduction}
Differential inclusions are a powerful instrument for modeling dynamical systems. In this paper the differential inclusions of some special structure are explored. More specific: the support function of the set in the right-hand side of a differential inclusion contains the sum of the maximum and the minimum functions of the finite number of continuously differentiable (in phase coordinates) functions. So it is required to find a trajectory of such a differential inclusion which simultaneously would satisfy the boundary conditions and lie on the surface prescribed.
\bigskip
\bigskip
\bigskip
\bigskip

Such differential inclusions arise from systems of differential equations with discontinuous right-hand sides (when considering such systems moving in a sliding mode), from linear control systems when the control function is subject to some nonsmooth (in phase vector) constraints or/and the nature of disturbance of the right-hand side is known to contain nondifferentiable functions of phase coordinates, some real mechanical models when coordinates, their velocities and accelarations are under some geometrical constraints leading to the differential inclusions structure considered. A good real example generating the type of differential inclusion under consideration is physical model of a system with dry friction and nonsmooth control forcing this system to move in a sliding mode (see Example 7.3 below). 


The majority of papers in literature considering differential inclusions are devoted to such classical problems as: existence of solutions
 \cite{lasop}, \cite{kikuchi}, \cite{Cellina}, \cite{smirnov}, dependence of solutions on the parameters \cite{smirnov},
\cite{clarke}, attainability and viability \cite{Cellina}, \cite{clarke}, \cite{krasquin}. Let us also give some references \cite {Cernea}, \cite{Pappas}, \cite{Zhu}, \cite{Aru} with optimality conditions in problems with differential inclusion. In these papers differential inclusions of a rather general form are considered. There are cases of phase constraints as well as nonsmooth and nonconvex ones. On the other hand, the works listed are more of theoretical significance and some results seem hard to be employed in practice. Note that the majority of methods in literature consider only differential inclusions with a free right endpoint and use some classical approaches as Euler and Runge-Kutta schemes, various finite differences methods etc. (see, e.~ g., \cite {Sandberg}, \cite{Bastien}, \cite{Beyn}, \cite{Lempio}, \cite{Veliov}). A survey of difference methods for differential inclusions can be found in \cite{Dontchev}. Note one paper \cite{Schilling} where an algorithm to solve boundary value problems for differential inclusions was constructed. 


We reduce the original problem to a variational one; such a reduction is implemented just like in the author previous works \cite{Fominyh_1},  \cite{Fominyh_2}. The functional obtained is proved to be quasidifferentiable. So the nonsmooth optimization methods are to be applied. In papers \cite{Fom_IJC}, \cite{Fom_opt_lett} a method for minimizing of a quasi (sub) differential functional is proposed based on the idea of considering phase trajectory and its derivative as independent variables (and taking the natural connection between these variables into account via the special penalty term). Also note the work \cite{Fom_NFA} where the similar technique is used to consider systems moving in sliding modes in a more simple case of smooth functionals. 

The ``weak'' convergence of the modified quasidifferential descent method is proved in some particlular cases due to the special structure of the functional to be minimized. The conceptual scheme of the method as well as ideas of proof are based on the analogous ones of V. F. Demyanov scientific school on nondifferentiable minimization and taken from book \cite{demmal} where a modified subdifferentiable (steepest) descent method for minimizing the maximum of the finite number of continuously differentiable functions in the finite dimensional space is proposed justified.

The principial novelity of the method suggested in the paper consists of the separation of the variables (phase coordinates and its derivatives) and in the subsequent results such as necessary minimum conditions in new effective form and construction of the quasidifferential descent method arising from this idea.

As it was noted above, only a small part of the literature, 
devoted to differential inclusions, deals with constructing numerical methods for solving the corresponding boundary value problem.
Moreover, to the best of the author knowledge literature devoted to differential inclusions all the more do not deal with the 
numerical methods for nonsmooth right-hand side of the inclusion as it is done in the current paper. 
 
\section{Basic definitions and notations} \label{sc2}
In this paper we use the following notations. Denote $\mathcal N$ the set of natural numbers. Let $C_{n} [0, T]$ be the space of $n$-dimensional continuous on $[0, T]$ vector-functions. Let also $P_{n} [0, T]$ be the space of piecewise continuous and bounded on $[0, T]$ $n$-dimensional vector-functions. We also require the space $L^2_n [0, T]$ of measurable and square-summable on $[0, T]$ $n$-dimensional vector-functions (strictly speaking the known corresponding factor-space is considered). If $X$ is some normed space, then $||\cdot||_X$ denotes its norm and $X^*$ denotes the space conjugate to the space given.
\newpage

We will assume that each trajectory $x(t)$ is a piecewise continuously differentiable vector-function. Let $t_0 \in [0, T)$ be a point of nondifferentiability of the vector-function $x(t)$, then we assume that $\dot x(t_0)$ is a right-hand derivative of the vector-function $x(t)$ at the point~$t_0$ for definiteness. Similarly, we assume that $\dot x(T)$ is a left-hand derivative of the vector-function $x(t)$ at the point~$T$. Now with the assumptions and the notations taken we can suppose that the vector function $x(t)$ belongs to the space $C_{n} [0, T]$ and that the vector function~$\dot x(t)$ belongs to the space $P_{n} [0, T]$.

For some arbitrary set $F \subset R^n$ define the support function of the vector $\psi~\in~R^n$ as $c(F, \psi) = \sup \limits_{f\in F}\langle f, \psi \rangle$ where $\langle a, b \rangle$ is a scalar product of the vectors $a, b \in R^n$. Denote $S_n$ and $B_n$ a unit sphere and a unit ball in $R^n$ with the center in the origin respectively, also let $B_r(c)$ be a ball with the radius $r \in R$ and the center $c \in R^n$. Let the vectors $\bf{e_i}$, $i =\overline{1,n}$, form the standard basis in $R^n$. Let $0_n$ denote a zero element of a functional space of some $n$-dimensional vector-functions and $\bf{0_n}$ denote a zero element of the space $R^n$.
If $\varphi_1(x) = \max\limits_{i = \overline{1, M}} f_i(x)$, where $f_i(x): R^n \rightarrow R$, $i = \overline{1, M}$, are some functions, then we call the function $f_{\overline i}(x)$, $\overline i \in \{1,\dots,M\}$, an active one at the point $x_0 \in R^n$, if $\overline i \in \mathcal R_1(x_0) = \Big\{i \in \{1,\dots,M\} \ \big| \ f_i(x_0) = \varphi_1(x_0) \Big\} $. If $\varphi_2(x) = \min\limits_{j = \overline{1, K}} g_j(x)$, where $g_j(x): R^n \rightarrow R$, $j = \overline{1, K}$, are some functions, then we call the function $g_{\overline j}(x)$, $\overline j \in \{1,\dots,K\}$, an active one at the point $x_0 \in R^n$, if $\overline j \in \mathcal R_2(x_0) = \Big\{j \in \{1,\dots,K\} \ \big| \ g_j(x_0) = \varphi_2(x_0) \Big\} $.


In the paper we will use both subdifferentials and superdifferentials of functions in a finite-dimensional space and subdifferentials and superdifferentials of functionals in a functional space. Despite the fact that the second concepts generalize the first ones, for convenience we separately introduce definitions for both of these cases and for those specific functions (functionals) and their variables and spaces which are considered in the paper. 

Consider the space $R^n \times R^n$ with the standard norm. Let $d = (d_1, d_2)' \in R^n \times R^n$ be an arbitrary vector. Suppose that at the point $(x, z)$ there exists such a convex compact set $\overline \partial h_1(x,z)$ $\subset R^n \times R^n$ that 
\begin{equation}
\label{0.1}
\frac{\partial h_1(x,z)}{\partial d} = \lim_{\alpha \downarrow 0} \frac{1}{\alpha} \Big(h_1(x+\alpha d_1, z + \alpha d_2) - h_1(x,z)\Big) = \min_{w \in \overline \partial h_1(x,z)} \langle w, d \rangle. 
\end{equation}

In this case the function $h_1(x,z)$ is called \cite{demmal} superdifferentiable at the point $(x, z)$, and the set $\overline \partial h_1(x,z)$ is called the superdifferential of the function $h_1(x,z)$ at the point $(x,z)$.

From expression (\ref{0.1}) one can see that the following formula holds true:
$$ h_1(x + \alpha d_1, z + \alpha d_2) = h_1(x, z) + \alpha  \min_{w \in \overline \partial h_1(x,z)} \langle w, d \rangle + o(\alpha, x, z, d), $$
$$\quad  \frac{o(\alpha, x, z, d)}{\alpha} \rightarrow 0, \ \alpha \downarrow 0.$$

Consider the space $R^n \times R^n$ with the standard norm. Let $d = (d_1, d_2)' \in R^n \times R^n$ be an arbitrary vector. Suppose that at the point $(x, z)$ there exists such a convex compact set $\underline \partial h_2(x,z)$ $\subset R^n \times R^n$ that 
\begin{equation}
\label{0.2}
\frac{\partial h_2(x,z)}{\partial d} = \lim_{\alpha \downarrow 0} \frac{1}{\alpha} \Big(h_2(x+\alpha d_1, z + \alpha d_2) - h_2(x,z)\Big) = \max_{v \in \underline \partial h_2(x,z)} \langle v, d \rangle. 
\end{equation}

In this case the function $h_2(x,z)$ is called \cite{demmal} subdifferentiable at the point $(x, z)$, and the set $\underline \partial h_2(x,z)$ is called the subdifferential of the function $h_2(x,z)$ at the point $(x,z)$.

From expression (\ref{0.2}) one can see that the following formula holds true:
$$ h_2(x + \alpha d_1, z + \alpha d_2) = h_2(x, z) + \alpha  \max_{v \in \underline \partial h_2(x,z)} \langle v, d \rangle + o(\alpha, x, z, d), $$
$$\quad  \frac{o(\alpha, x, z, d)}{\alpha} \rightarrow 0, \ \alpha \downarrow 0.$$

%

If the function $\varsigma(\xi)$ is differentiable at the point $\xi_0 \in R^\ell$, then its superdifferential (subdifferential) at this point is represented \cite{demmal} in the form \begin{equation}
\label{0.4}\overline \partial \varsigma(\xi_0) = \{ \varsigma'(\xi_0) \}, \quad \underline \partial \varsigma(\xi_0) = \{ \bf{0_\ell} \}
\end{equation}
$$
\Big( \underline \partial \varsigma(\xi_0) = \{ \varsigma'(\xi_0) \}, \quad \overline \partial \varsigma(\xi_0) = \{ \bf{0_\ell} \} \Big),
$$
 where $\varsigma'(\xi_0)$ is a gradient of the function $\varsigma(\xi)$ at the point $\xi_0$.
 
Note also that the superdifferential (subdifferential) of the finite sum of superdifferentiable (subdifferentiable) functions is the sum of the superdifferentials (subdifferentials) of summands, i. e. if the functions $\varsigma_k(\xi)$, $k = \overline{1, r}$, are superdifferentiable (subdifferentiable) at the point $\xi_0 \in R^\ell$, then the function $\varsigma(\xi) = \sum_{k=1}^r \varsigma_k(\xi)$ superdifferential (subdifferential) at this point is calculated \cite{demmal} by the formula 
\begin{equation}
\label{0.5} \overline \partial \varsigma(\xi_0) = \sum_{k=1}^r \overline\partial \varsigma_k(\xi_0).
\end{equation}
$$
\Big( \underline \partial \varsigma(\xi_0) = \sum_{k=1}^r \underline\partial \varsigma_k(\xi_0) \Big).
$$

Consider the space $C_n[0, T] \times P_n[0, T]$ with the norm $L_n^2 [0, T] \times L_n^2 [0, T]$. Let $g = (g_1, g_2)' \in C_n[0, T] \times P_n[0, T]$ be an arbitrary vector-function. Suppose that at the point $(x, z)$ there exists such a convex weakly$^*$ compact set $\overline \partial {\varphi_1(x, z)} \subset \big( C_n[0, T] \times P_n[0, T],$ $|| \cdot ||_{L_n^2 [0, T] \times L_n^2 [0, T]} \big) ^*$ that  
\begin{equation}
\label{0.6} 
\frac{\partial \varphi_1(x, z)}{\partial g} = \lim_{\alpha \downarrow 0} \frac{1}{\alpha} \Big(\varphi_1(x+\alpha g_1, z+\alpha g_2) - \varphi_1(x, z)\Big) = \min_{w \in \overline \partial \varphi_1(x, z)} w(g).
\end{equation}

In this case the functional $\varphi_1(x,z)$ is called \cite{dolgquasi} superdifferentiable at the point $(x, z)$ and the set $\overline \partial {\varphi_1(x, z)} $ 
is called a superdifferential of the functional $\varphi_1(x, z)$ at the point $(x,z)$. 

From expression (\ref{0.6}) one can see that the following formula holds true:
$$ 
\varphi_1(x + \alpha g_1, z + \alpha g_2) = \varphi_1(x, z) + \alpha  \min_{w \in \overline \partial \varphi_1(x, z)} w(g) + o(\alpha, x, z, g),
$$
$$ \quad  \frac{o(\alpha, x, z, g)}{\alpha} \rightarrow 0, \ \alpha \downarrow 0.$$

Consider the space $C_n[0, T] \times P_n[0, T]$ with the norm $L_n^2 [0, T] \times L_n^2 [0, T]$. Let $g = (g_1, g_2)' \in C_n[0, T] \times P_n[0, T]$ be an arbitrary vector-function. Suppose that at the point $(x, z)$ there exists such a convex weakly$^*$ compact set $\underline \partial {\varphi_2(x, z)} \subset \big( C_n[0, T] \times P_n[0, T],$ $|| \cdot ||_{L_n^2 [0, T] \times L_n^2 [0, T]} \big) ^*$ that  
\begin{equation}
\label{0.7} 
\frac{\partial \varphi_2(x, z)}{\partial g} = \lim_{\alpha \downarrow 0} \frac{1}{\alpha} \Big(\varphi_2(x+\alpha g_1, z+\alpha g_2) - \varphi_2(x, z)\Big) = \max_{v \in \underline \partial \varphi_2(x, z)} v(g).
\end{equation}

In this case the functional $\varphi_2(x,z)$ is called \cite{dolgquasi} subdifferentiable at the point $(x, z)$ and the set $\underline \partial {\varphi_2(x, z)} $ 
is called a subdifferential of the functional $\varphi_2(x, z)$ at the point $(x,z)$. 

From expression (\ref{0.7}) one can see that the following formula holds true: 
$$ 
\varphi_2(x + \alpha g_1, z + \alpha g_2) = \varphi_2(x, z) + \alpha  \max_{v \in \underline \partial \varphi_2(x, z)} v(g) + o(\alpha, x, z, g),
$$
$$ \quad  \frac{o(\alpha, x, z, g)}{\alpha} \rightarrow 0, \ \alpha \downarrow 0.$$

\begin{rmrk}
Note that more general definitions (both for the finite-dimensional and the functional spaces) may be given, when the function (functional) directional derivative may be represented as a sum of a maximum and a minimum of linear functionals over some sets (the subdifferential and the superdifferential respectively), then the function (functional) is called quasidifferentiable and the pair [subdifferential, superdifferential] is called a quasidifferential. However we intentionally give the independent definitions of these objects, since the support function under consideration in this paper is a sum of a superdifferentiable and a subdifferentiable functions (see formula (\ref{support}) below).  
\end{rmrk}

\section{Statement of the problem and reduction to a variational one} 
 We now give the formal statement of the problem. The motivating examples leading to such a problem are given in Section 5. 

It is required to find such a trajectory $x^{*} \in C_{n}[0, T]$ (with the derivative $\dot x^{*} \in P_{n}[0, T]$) satisfies differential inclusion
\begin{equation} \label{4}
\dot x(t) \in F(x(t)),
\end{equation}
moves along surface
\begin{equation}
\label{5}
e(x) = {\bf 0_o}
\end{equation}
 while $t \in [0,T]$ and meets boundary conditions
\begin{equation}
\label{2}
x(0) = x_{0},
\end{equation}
\begin{equation}
\label{3}
x_j(T) = {x_T}_j, \quad j \in J.
\end{equation}
Assume that there exists such a solution.

In (\ref{5}) let $e(x)$ be a known continuously differentiable $o$-dimensional vector-function; in formula (\ref{2}) the initial point $x_0 \in R^n$ is a given vector; in formula (\ref{3}) the final point coordinates ${x_T}_j$, $j \in J$, are given numbers corresponding to those ones of the state vector which are fixed at the right endpoint, here $J \subset \{1, \dots, n\}$ is a given index set.

Now make assumptions regarding the right-hand side support function. Suppose that differential inclusion in~(\ref{4}) is endowed with the following properties: 1) at each fixed $x \in R^n$ the set $F(x)$ on the right-hand side may be ``partitioned'' into the sets $F_1(x), \dots F_n(x)$ (such that $\dot x_i \in F_i(x)$, $i = \overline{1,n}$) and each of the sets $F_1(x), \dots F_n(x)$ is convex and compact; 2) the support functions of the corresponding sets $F_i(x)$, $i = \overline{1,n}$, may be represented as 
\begin{equation} \label{support} c(F_i(x),\psi_i) = c_1(F_i(x),\psi_i) + c_2(F_i(x),\psi_i), \end{equation}
where
\begin{equation} \label{supf}  c_1(F_i(x),\psi_i) = \sum_{j=1}^{r} \max\Big\{f_{i, j_1}(x) \psi_i, \dots, f_{i, j_{m(j)}}(x) \psi_i\Big\}, \end{equation}
\begin{equation} \label{subf} c_2(F_i(x),\psi_i) = \sum_{j=1}^{s} \min\Big\{g_{i, j_1}(x,  \psi_i), \dots, g_{i, j_{k(j)}}(x, \psi_i) \Big\}, \end{equation}
where $f_{i, j_1}(x) \psi_i, \dots,  f_{i, j_{m(j)}}(x) \psi_i$, $i = \overline{1, n}$, $j = \overline{1, r}$, $g_{i, j_1}(x, \psi_i), \dots,  g_{i, j_{k(j)}}(x, \psi_i)$, $i = \overline{1, n}$, $j = \overline{1, s}$ (for simplicity of presentation we
suppose that $r$ and $s$ are the same for each $i = \overline{1, n}$), are continuously differentiable in $x$ functions (at fixed $\psi_i \in S_1$, $i = \overline{1,n}$). 

A slight difference in the form of formulas of two support functions above is due to practical problems generating such a structure of support functions (see Section 5 for corresponding examples).

\bigskip
We will sometimes write $F$ instead of $F(x)$ for brevity. 
Insofar as  $\forall x \in R^n$ the set $F_i(x)$, $i = \overline{1,n}$, 
is a convex compact set in $R^n$, then inclusion (\ref{4}) 
may be rewritten {\cite{Blagodatskih}} as follows:
$$
\dot x_i(t) \psi_i(t) \leq c(F_i(x(t)), \psi_i(t)) \quad \forall \psi_i(t) \in S_1, \quad \forall t \in [0, T], \quad i = \overline{1,n}.
$$

Denote $z(t) = \dot x(t)$, $z \in P_n [0, T]$, then from (\ref{2}) one has
\begin{equation}
\label{3.11''}
x(t) = x_0 + \int_0^t z(\tau) d\tau.
\end{equation}

Let us now realize the following idea. ``Forcibly'' consider the points $z$ and $x$ to be ``independent'' variables. Since, in fact, there is relationship (\ref{3.11''}) between these variables (which naturally means that the vector-function~$z(t)$ is a derivative of the vector-function $x(t)$), let us take this restriction into account by using the functional
 \begin{equation} 
\label{3.17}
 \upsilon(x, z) = \frac{1}{2} \int_0^T \left( x(t) - x_0 - \int_0^t z(\tau) d \tau \right)^2 dt.
\end{equation} 
It is seen that relation (\ref{3.11''}) and condition (\ref{2}) on the left endpoint is satisfied iff $\upsilon(x, z) = 0$. 

For $i = \overline{1,n}$ put
$$
\ell_i(\psi_i, x, z) = \langle z_i, \psi_i \rangle - c ( F_i(x), \psi_i ),
$$
$$
h_i(x, z) = \max_{\psi_i \in S_1} \max \{ 0, \ell_i(\psi_i, x, z) \},
$$
then put
$$h(x,z) = (h_1(x,z), \dots, h_n(x,z))'$$
and construct the functional
\begin{equation}
\label{3.13}
\varphi(x, z) = \frac{1}{2} \int_0^T || h \big( x(t), z(t) \big) ||^2_{R^n} dt.
\end{equation}

It is not difficult to check that for functional (\ref{3.13}) the relation
$$
\left\{
\begin{array}{ll}
\varphi(x, z) = 0, \ &\text{if} \ \dot x_i(t) \psi_i(t)  \leq c(F_i(x(t)), \psi_i(t)) \quad \forall \psi_i(t) \in S_1, \quad \forall t \in [0, T], \quad i = \overline{1,n}. \\
\varphi(x, z) > 0, \ &\text{otherwise},
\end{array}
\right.
$$
holds true, i. e. inclusion (\ref{4}) takes place iff $\varphi(x, z) = 0$.

Introduce the functional
\begin{equation}
\label{3.14}
\chi(z) = \frac{1}{2} \sum_{j \in J} \left( {x_0}_j + \int_0^T z_j(t) dt - {x_T}_j \right)^2.
\end{equation}

 We see that if $\upsilon (x, z) = 0$, then condition (\ref{3}) on the right endpoint is satisfied iff $\chi(z) = 0$.

Introduce the functional 
\begin{equation}
\label{3.14'}
\omega(x) = \frac{1}{2} \int_0^T ||e \big( x(t) \big )||^2_{R^o} dt.
\end{equation}

Obviously, the trajectory $x(t)$ belongs to surface (\ref{5}) at each $t \in [0, T]$ iff $\omega(x) = 0$.

Finally construct the functional
\begin{equation} 
\label{3.16} 
I(x, z) = \varphi(x, z) + \chi(z) + \omega(x) + \upsilon(x, z).
\end{equation}

So the original problem has been reduced to minimizing functional (\ref{3.16}) on the space $$X = \Big( C_n[0, T] \times P_n [0, T],  L_n^2 [0, T] \times L_n^2 [0, T] \Big).$$ Denote $z^*$ a global minimizer of this functional. Then $$ x^*(t) = x_0 + \int_0^t z^*(\tau) d\tau $$
is a solution of the initial problem. 

\begin{rmrk} \label{rm4.1}
The structure of the functional $\varphi(x, z)$ is natural as the value $ h_i(x(t), z(t)) $, $i = \overline{1, n}$,
at each fixed $t \in [0, T]$ is the Euclidean distance from the point $z_i(t)$
to the set $F_i (x(t))$; functional (\ref{3.13}) is half the sum of squares of the deviations in $L^2_n[0,T]$ norm of the trajectories $z_i(t)$ from the sets $F_i(x)$, $i = \overline{1, n}$, respectively. The meaning of functionals (\ref{3.17}), (\ref{3.14}), (\ref{3.14'}) structures is obvious.
\end{rmrk}

\begin{rmrk} \label{rm4.2}
Despite the fact that the dimension of functional $I(x, z)$ arguments is $n$ more the dimension of the initial problem (i. e. the dimension of the point $x^*$), the structure of its superdifferential (in the space $C_n[0, T] \times P_n[0, T]$ as a normed space with the norm $L_n^2 [0, T] \times L_n^2 [0, T]$), as will be seen from what follows, has a rather simple form. This fact will allow us to construct
a numerical method for solving the original problem. 
\end{rmrk}

\section{Minimum conditions of the functional $I(x, z)$}
In this section referring to superdifferential calculus rules (\ref{0.4}), (\ref{0.5}), we mean their known analogues in a functional space \cite{Dolgopolikcodiff}. 

Using classical variation it is easy to prove the Gateaux differentiability of the functional $\chi(z)$, we have
$$
 \nabla \chi(z) = \sum_{j \in J} \left( {x_0}_j + \int_0^T z_j(t) dt - {x_T}_j \right) {\bf e_j}.  
 $$
 By superdifferential calculus rule (\ref{0.4}) one may put 
\begin{equation} \label{4.44} \underline \partial \chi(z) =  \left\{ \sum_{j \in J} \left( {x_0}_j + \int_0^T z_j(t) dt - {x_T}_j \right) {\bf e_j} \right\}, \quad \overline \partial \chi(z) = 0_n \end{equation}
 
Using classical variation it is easy to prove the Gateaux differentiability of the functional $\omega(x)$, we have
$$
 \nabla \omega(x, t) = \sum_{j=1}^o e_j\big(x(t)\big) \frac{\partial e_j\big(x(t)\big)} {\partial x}.  
 $$
 By superdifferential calculus rule (\ref{0.4}) one may put 
\begin{equation} \label{4.6} \underline \partial \omega(x, t)  =  \left\{ \sum_{j=1}^o e_j\big(x(t)\big) \frac{\partial e_j\big(x(t)\big)} {\partial x} \right\}, \quad \overline \partial \omega(x, t)  = 0_n. \end{equation}

 Using classical variation and integration by parts it is also not difficult to check the Gateaux differentiability of the functional $\upsilon(x, z)$, we obtain
 $$\nabla \upsilon(x, z, t) = \begin{pmatrix}
\displaystyle{ x(t) - x_0 - \int_0^t z(\tau) d\tau } \\
\displaystyle{ -\int_t^T \left( x(\tau) - x_0 - \int_0^\tau z(s) ds \right) d\tau  }
\end{pmatrix}. 
$$   
 
  By superdifferential calculus rule (\ref{0.4}) one may put 
\begin{equation} \label{4.55} \underline \partial \upsilon(x, z, t)  =  \left\{ \begin{pmatrix}
\displaystyle{ x(t) - x_0 - \int_0^t z(\tau) d\tau } \\
\displaystyle{ -\int_t^T \left( x(\tau) - x_0 - \int_0^\tau z(s) ds \right) d\tau  }
\end{pmatrix} \right\}, \quad \overline \partial \upsilon(x, z, t) = 0_{2n}. \end{equation}

\begin{rmrk}
Note the following fact. Since, as is known, the space $\Big( C_n[0, T], || \cdot ||_{L^2_n [0, T]} \Big)$ is everywhere dense in the space ${L^2_n [0, T]}$ and the space $\Big( P_n[0, T], || \cdot ||_{L^2_n [0, T]} \Big)$ is also everywhere dense in the space ${L^2_n [0, T]}$, then the space $X^*$ conjugate to the space $X$ introduced in the previous paragraph is isometrically isomorphic to the space ${L^2_n [0, T] \times L^2_n [0, T]}$ (see \cite{kolfom}). Now keeping this fact in mind, note that strictly speaking, in formulas (\ref{4.44}), (\ref{4.6}), (\ref{4.55}) not the subdifferentials and superdifferentials theirselves but their images under the corresponding isometric isomorphisms are given.
\end{rmrk}

Explore the differential properties of the functional $\varphi(x,z)$. For this, we first give the following formulas for calculating the quasidifferential $\Big[ \underline \partial h^2(x, z), \overline \partial h^2(x, z) \Big]$ at the point $(x, z)$. At $i = \overline{1,n}$ one has
\begin{equation}
\label{4.5}
\overline \partial \Big( {\textstyle \frac{1}{2}} \, h_i^2(x, z) \Big) = h_i(x,z)\bigg(\psi^*_i {\bf e_{i+n}} + \sum_{j=1}^{r} \overline \partial \Big(-\max\Big\{f_{i, j_1}(x) \psi_i^*, \dots, f_{i, j_{m(j)}}(x) \psi_i^*\Big\} \Big)\bigg)
\end{equation}

where at $j = \overline{1,r}$ and $j_q \in \mathcal R^1_{ij}(x)$ we have
$$ \overline \partial \Big(-\max\Big\{f_{i, j_1}(x) \psi_i^*, \dots, f_{i, j_{m(j)}}(x)   \psi_i^*\Big\} \Big) = \mathrm{co} \left \{ 
\left[ -\psi_i^* \frac{\partial f_{i, j_q}(x)}{\partial x}, \bf{0_n} \right] 
\right \},
$$
$$
\mathcal R^1_{ij}(x) = \Big\{j_q \in \{j_1, \dots, j_{m(j)}\} \ \big| \ f_{i, j_q}(x) \psi_i^* =  \max\Big\{f_{i, j_1}(x) \psi_i^*, \dots, f_{i, j_{m(j)}}(x) \psi_i^*\Big\} \Big\}
$$
and
\begin{equation}
\label{4.555}
\underline \partial \Big( {\textstyle \frac{1}{2}} \, h_i^2(x, z) \Big) = h_i(x,z) \sum_{j=1}^{s} \underline \partial \Big(-\min\Big\{g_{i, j_1}(x, \psi_i^*), \dots, g_{i, j_{k(j)}}(x, \psi_i^*) \Big\} \Big)
\end{equation}
where at $j = \overline{1,s}$ and  $j_p \in \mathcal R^2_{ij}(x)$ we have
$$
\underline \partial \Big(- \min\Big\{g_{i, j_1}(x, \psi_i^*), \dots, g_{i, j_{k(j)}}(x, \psi_i^*) \Big\} \Big) =  \ \mathrm{co} \left \{ 
\left[ -\frac{\partial g_{i, j_p}(x, \psi_i^*)}{\partial x}, \bf{0_n} \right] 
\right \},
$$
$$
\mathcal R^2_{ij}(x) = \bigg\{j_p \in \{j_1, \dots, j_{k(j)}\} \ \big| \ g_{i, j_p}(x, \psi_i^*) =  \min\Big\{g_{i, j_1}(x, \psi_i^*), \dots, g_{i, j_{k(j)}}(x, \psi_i^*) \Big\} \bigg\}
$$
and
$$
\ell_i(\psi_i^*(x, z), x, z) =  \max_{\psi_i \in S_1} \ell_i(\psi_i, x, z) \quad if \  h_i(x, z) > 0 $$
$$ \psi_i^* = \psi_i^0 \ ({where} \ \psi_i^0 \ is \ some \ fixed \ element \ of \ the \ set \ S_1) \quad if \ h_i(x, z) = 0.$$

Note that in the case $h_i(x, z) > 0$ the element $\psi_i^*(x, z)$ is unique from support functions properties \cite{Blagodatskih} and is continuous \cite{BonnansShapiro} from maximum (minimum) function properties.

Formulas (\ref{4.5}) and (\ref{4.555}) can be easily checked based on the proof given in \cite{Fom_NFA} and on the quasidifferential calculus rules \cite{demvas}. Therefore, let us give only a brief scheme of proof here.

 Fix some $i \in \{1, \dots, n\}$ and consider two cases. 

a) Suppose that $h_i (x, z) > 0$, i. e. $h_i(x, z) = \max_{\psi_i \in S_1}\ell_i (\psi_i, x, z) > 0$. 
Our aim is to apply the corresponding theorem on a directional differentiability from \cite{BonnansShapiro}. There the inf-functions are considered, hence let us apply this theorem to the function $-\ell_i(\psi_i, x, z)$. Check that the function $h_i(x, z)$ satisfies the conditions: \newline
i) the function $\ell_i (\psi_1, x, z)$ is continuous on $S_1 \times R^n \times R^n$; \newline
ii)  there exist a number $\beta$ and a compact set $C \in R$ such that for every $(\overline x, \overline z)$ in the vicinity of the point $(x, z)$ the level set 
${lev}_\beta (-\ell(\cdot, \overline x, \overline z)) = \{\psi_i \in S_1 \ | \ -\ell_i(\psi_i, \overline x, \overline z) \leq \beta\} $
is nonempty and is contained in the set $C$; \newpage
iii) for any fixed $\psi_i \in S_1$ the function $\ell_i(\psi_i, \cdot, \cdot)$ is directionally differentiable at the point $(x, z)$; \newline
iv) if $d = (d_1, d_2)' \in R^n \times R^n$, $\gamma_n \downarrow 0$ and $\psi_{i_{n}}$ is a sequence in $C$, then $\psi_{i_{n}}$ has a limit point $\overline \psi_i$ such that 
$$\mathrm{lim} \sup\limits_{n \rightarrow \infty} \frac{-\ell_i(\psi_{i_{n}}, x + \gamma_n d_{1},  z+\gamma_n d_2) - (-\ell_i(\psi_{i_{n}}, x, z))}{\gamma_n} \geq \frac{\partial(-\ell(\overline \psi_i, x, z))}{\partial d}. $$

The verification of conditions i), ii) is obvious. 

In order to verify iii) it is sufficient to note that at every fixed $\psi_i \in S_1$ the corresponding functions under consideration are the maxima (minima) of continuously differentiable functions and it is well known \cite{demvas} that the $\max$-function ($\min$-function) (of continuously differentiable functions) is subdifferentiable (superdifferentiable), hence the explicit expression for their directional derivatives can be obtained.

In order to verify iv) it is sufficient to write out explicitly the right-hand and the left-hand sides of the inequality required and to check directly that iv) holds true using the Taylor expansion of the continuously differentiable functions under $\max$-function ($\min$-function) and the well known \cite{demmal} inequalities $$\min_{j = \overline{1,K}}(\nu_j + \mu_j) \geq \min \nu_{j = \overline{1,K}} + \min \mu_{j = \overline{1,K}},$$ $$\max_{i = \overline{1,M}}(\nu_i + \mu_i) \geq \max_{i = \overline{1,M}} \nu_i + \max_{i \in \mathcal R} \mu_i, \quad \mathcal R = \Big\{i \in \{1, \dots, M\} \ | \ \nu_i = \max_{i = \overline{1,M}} \nu_i \Big\},$$ which are true for the arbitrary numbers $\nu_i$, $\mu_i$, $i = \overline{1, K}$ or $\nu_i$, $\mu_i$, $i = \overline{1, M}$.

So the function $h_i(x, z)$ is directionally differentiable at the point $(x, z)$ but writing out explicitly its directional derivative we make sure that by definition (see (\ref{0.1})) the function $h_i(x, z)$ is quasidifferentiable at the point $(x, z)$; so the function $h_i^2(x, z)$ is quasidifferentiable at the point $(x, z)$ and its subdifferential and superdifferential are obtained via (\ref{4.555}) and (\ref{4.5}) respectively by quasidifferential calculus rules \cite{demvas}.

b) In the case $h_i (x, z) = 0$ it is obvious that the function $h_i^2(x, z)$ is differentiable at the point $(x, z)$, and its gradient vanishes at this point.  


If the interval $[0, T]$ may be divided into a finite number of intervals, in every of which one (several) of the considered functions is (are) active, then the function $h_i(x, z)$, $i = \overline{1,n}$, quasidifferentiability can be proved on every of such intervals in a completely analogous fashion. In \cite{Fom_IJC} it is shown that the functional $\varphi(x,z)$ is quasidifferentiable and that its quasidifferential is determined by the corresponding integrand quasidifferential. 

\begin{thrm}  \label{th4.2}
Let the interval $[0, T]$ be divided into a finite number of intervals, in every of which one (several) of the functions $\big\{f_{i, j_1}(x) \psi_i, \dots, f_{i, j_{m(j)}}(x) \psi_i\big\}$, $i = \overline{1,n}$, $j = \overline{1,r}$, $\big\{g_{i, j_1}(x, \psi_i), \dots, g_{i, j_{k(j)}}(x, \psi_i)\big\}$, $i = \overline{1,n}$, $j = \overline{1,s}$, is (are) active. Then the functional
 $
 \varphi(x, z) 
 $
is quasidifferentiable, i. e. 
\begin{equation}
\label{11.10} 
\frac{\partial \varphi(x, z)}{\partial g} = \lim_{\alpha \downarrow 0} \frac{1}{\alpha} \Big(\varphi(x+\alpha g_1, z+\alpha g_2) - \varphi(x, z)\Big) = \max_{v \in \underline \partial \varphi(x, z)} v(g) + \min_{w \in \overline \partial \varphi(x, z)} w(g) \end{equation}
Here the sets $\underline \partial \varphi(x, z)$, $\overline \partial \varphi(x, z)$ are of the following form
\begin{equation}
\label{11.2001}  \underline \partial \varphi(x, z) = \Bigg\{ v \in \Big(C_n[0, T] \times P_n[0, T], || \cdot ||_{L^2_n[0,T] \times L^2_n[0,T]} \Big)^* \ \Big| \end{equation}  
$$v(g) = \int_0^T \langle v_1(t), g_1(t) \rangle dt + \int_0^T \langle v_2(t), g_2(t) \rangle dt \quad \forall g_1 \in C_n[0,T], \, g_2 \in P_n[0,T], $$
$$ v_1(t), v_2(t) \in L_n^\infty[0, T], \quad (v_1(t), v_2(t))' \in {\textstyle \frac{1}{2}} \, \underline \partial h^2(x(t), z(t)) \quad for \ all \ t \in [0, T] \Bigg\}.$$

\begin{equation}
\label{11.2000}  \overline \partial \varphi(x, z) = \Bigg\{ w \in \Big(C_n[0, T] \times P_n[0, T], || \cdot ||_{L^2_n[0,T] \times L^2_n[0,T]} \Big)^* \ \Big| \end{equation}  
$$w(g) = \int_0^T \langle w_1(t), g_1(t) \rangle dt + \int_0^T \langle w_2(t), g_2(t) \rangle dt \quad \forall g_1 \in C_n[0,T], \, g_2 \in P_n[0,T], $$
$$ w_1(t), w_2(t) \in L_n^\infty[0, T], \quad (w_1(t), w_2(t))' \in {\textstyle \frac{1}{2}} \, \overline \partial h^2(x(t), z(t)) \quad for \ all \ t \in [0, T] \Bigg\}.$$
\end{thrm}

Using formulas (\ref{4.44}), (\ref{4.6}), (\ref{4.55}), (\ref{4.5}), (\ref{4.555}), (\ref{11.10}), (\ref{11.2001}), (\ref{11.2000}) obtained and superdifferential calculus rule (\ref{0.5}) we have the final formula for calculating the subdifferential and the superdifferential of the functional $I(x, z)$ at the point $(x, z)$, let us denote
$$
\underline \partial \mathcal{I}(x, z, t) := {\textstyle \frac{1}{2}} \underline \partial h^2(x(t), z(t)) + \underline \partial \overline \chi(z) + \underline \partial \overline \omega(x, t) + \underline \partial \upsilon(x,z,t), \quad \overline \partial \mathcal{I}(x, z,t) := {\textstyle \frac{1}{2}} \overline \partial h^2(x(t), z(t)). $$

Using the known minimum condition {\cite{dolgquasi}} of the functional $I(x, z)$ at the point $(x^*, z^{*})$ in terms of quasidifferential, we conclude that the following theorem is true.

\begin{thrm}\label{th11.2}
Let the interval $[0, T]$ be divided into a finite number of intervals, in every of which one (several) of the functions $\big\{f_{i, j_1}(x) \psi_i, \dots, f_{i, j_{m(j)}}(x) \psi_i\big\}$, $i = \overline{1,n}$, $j = \overline{1,r}$, $\big\{g_{i, j_1}(x, \psi_i), \dots, g_{i, j_{k(j)}}(x,  \psi_i)\big\}$, $i = \overline{1,n}$, $j = \overline{1,s}$, is (are) active. In order for the point $(x^{*}, z^*)$ to minimize the functional $I(x, z)$, it is necessary to have 
\begin{equation}
\label{11.400}
-\overline \partial \mathcal{I} (x^*, z^{*}, t) \subset \underline \partial \mathcal{I} (x^*, z^{*}, t)
\end{equation}
at almost each $t \in [0,T]$. 

If one has $I(x^*, z^{*}) = 0$, then condition (\ref{11.400}) is also sufficient.
\end{thrm}


Theorem \ref{th11.2} already contains a constructive minimum condition, since on its basis it is possible to construct the quasidifferential descent direction, and for solving each of the subproblems arising during this construction there exist some known efficient algorithms for solving them.

\section{Problems Generating Nonsmooth Differential Inclusions} \label{sc09}

Consider the differntial inclusion
$$
\dot{x}(t) \in F(x(t)). 
$$

Let us consider the case when the right-hand side of the differential inclusion is the convex hull of several sets. For simplicity of presentation, we will describe the case with only two sets. Consider the mapping $$F(x) = \mathrm{co} \{F_1(x), F_2(x) \},$$ where the multivalued mappings $F_1(x), F_2(x)$ are continuous and at each fixed $x \in R^n$ the sets $F_1(x), F_2(x)$ are convex and compact. Using a simple formula \cite{Blagodatskih} for the support function of the union of two convex compact sets, we have $$c(F(x), \psi) = \max\Big\{c(F_1(x), \psi), c(F_2(x), \psi)\Big\}.$$  

Give here one one practical problem which leads to such a system. Let from some physical considerations the ``velocity'' $\dot x_1$ of an object lie in the range $[\min\{x_1, x_2, x_3\}, \max\{x_1, x_2, x_3\}]$ of the ``coordinates'' $x_1$, $x_2$, $x_3$. The segment given may be written down as $\mathrm{co} \{x_1, x_2, x_3\} = \mathrm{co}\bigcup\limits_{i=1}^3 \mathcal F_1^i(x)$, where $\mathcal F_1^i(x) = \{x_i\}$, $i = \overline{1,3}$. The support function \cite{Blagodatskih} of this set for $\psi_1 \in R$ is $$c(F(x), \psi) = \max \Big\{c\big(\mathcal F_1^1(x), \psi_1\big), c\big(\mathcal F_1^2(x), \psi_1\big), c\big(\mathcal F_1^3(x), \psi_1\big) \Big\} = \max  \{x_1 \psi_1, x_2 \psi_1, x_3 \psi_1\}.$$ 

Let us consider the case when the right-hand side of the differential inclusion is the convex hull of several sets. For simplicity of presentation, we will describe the case with only two sets. Consider the mapping $$F(x) = {F_1(x) \cap F_2(x) },$$ where the multivalued mappings $F_1(x), F_2(x)$ are continuous and at each fixed $x \in R^n$ the sets $F_1(x), F_2(x)$ are convex and compact. Using a known formula \cite{Blagodatskih} for the support function of the intersection of two convex compact sets, we have \begin{equation} \label{sf} c(F(x), \psi) = \overline{\mathrm {co}}\min\Big\{c(F_1(x), \psi), c(F_2(x), \psi)\Big\}. \end{equation}

Give here one practical problem where such differential inclusions may arise. Let from some physical considerations the ``velocity'' $\dot x_1$ of an object must belong both to the ranges $\mathcal F_1^1(x) = [x_1, x_2]$ and $\mathcal F_1^2(x) = [x_2-1, x_2+1]$ of the ``coordinates'' $x_1$, $x_2$. The set considered may be written down as $\mathcal F_1^1(x) \cap \mathcal F_1^2(x)$. The support function of this set for $\psi_1 \in S_1$ is \cite{Blagodatskih} as follows: $$ c(F(x), \psi) = \min \Big\{c\big(\mathcal F_1^1(x), \psi_1\big), c\big(\mathcal F_1^2(x), \psi_1\big)\Big\} = \min \left\{\frac{x_1+x_2}{2}\psi_1 + \frac{|x_1-x_2|}{2}, \, x_2 \psi_1 + 1\right\}, $$
and the first function under the minimum is continuously differentiable in $x$ if $x_1 \neq x_2$.

\begin{rmrk} \label{rm5.1} Note the following fact. As seen in example above in some cases one can directly check that for $\psi \in S_n$ the support function of the intersection of two sets is the minimum of the support functions of these sets. However, it is obvious that in general case the right-hand side of formula (\ref{sf}) does not coincide with $$\min\Big\{c(F_1(x), \psi), c(F_2(x), \psi)\Big\}$$ (there may exist $\psi \in R^n$ at which the values of these expressions are different). On the other hand using the known property of the support functions, it is easy to show that the point $z$ belongs to the intersection $F(x)$ iff $$\langle \psi, z \rangle \leq \min\Big\{c(F_1(x), \psi), c(F_2(x), \psi)\Big\} \quad \forall \psi \in S_n$$
  
This statement can be expanded to a more general case when the right-hand side is the sum of the finite number of the unions and intersections of the finite number of convex compacts. For simplicity we write out this condition for the sum of the union $F^1(x) = F_1^1(x) \cup F_2^1(x)$ and the intersection  $F^2(x) = F_1^2(x)\cap F_2^2(x)$: the point $z$ belongs to the set $F^1(x)+F^2(x)$ iff $$\langle \psi, z \rangle \leq \max\Big\{c(F_1^1(x),\psi), c(F_2^1(x), \psi)\Big\} + \min\Big\{c(F_1^2(x),\psi), c(F_2^2(x), \psi)\Big\} \quad \forall \psi \in S_n.$$

Hence in practice we may use this condition of belonging of a point to the set, we actually need not the support function of the intersection, but the support functions of the sets theirselves at $\psi \in S_n$. So instead of the functional $I(x,z)$ one may construct the similar functional  
$$\widetilde I(x, z) = \widetilde \varphi(x, z) + \chi(z) + \omega(x) + \upsilon(x, z),$$
where
$$
\widetilde \varphi(x, z) = \frac{1}{2} \int_0^T h^2 \big( x(t), z(t) \big) dt, \quad 
$$
$$
\widetilde \ell(\psi, x, z) = \langle z, \psi \rangle - \widetilde c ( F(x), \psi ), \quad 
\widetilde h(x, z) = \max_{\psi \in S_n} \max \{ 0, \widetilde\ell(\psi, x, z) \},
$$
and $\widetilde c ( F(x), \psi )$ is still the sum of maxima and minima of corresponding continuously differentiable in $x$ functions. 

However, when such a $\widetilde \varphi (x,z)$ functional structure is used, the following drawback arises: since $\widetilde c ( F(x), \psi )$, is no more a support function, the crucial property of the uniqueness of the vector $\psi^*$, is no more valid (in the case when the point $z$ does not belong to the set considered. Roughly speaking this is explained by the fact that function $\widetilde h(x,z)$, is no more the distance to the intersection of the sets but the maximum of the distances to these sets. While the distance to the convex compact (as the intersection of convex compacts) is achieved at the only point of the set, the maximum of the distances to these sets may be achieved at two points of different sets (if the distances are equal). Let us illustrate this by the following example. Take $F_1 = B_1((-1,0)')$, $F_2 = B_1((1,0)')$. Then one has $c(F_1, \psi) = -\psi_1 + 1$, $c(F_2, \psi) = \psi_1 + 1$. If one takes $z = (0, 1)'$, then $\max_{\psi \in S_2} \Big\{\psi_2 - \min\{-\psi_1 + 1, \psi_1 + 1\} \Big\}$ (which is the distance equal to $\sqrt{2} - 1$ of the vector $z$ both to the sets $F_1$ and $F_2$) is achieved at two elements: $\psi^* = (\sqrt{2}/2, \sqrt{2}/2)'$ or $\psi^* = (-\sqrt{2}/2, \sqrt{2}/2)'$.  

Hence in order for the Theorem \ref{th4.2} to remain true, we must make an additional assumption. We suppose that in the case $z \notin F(x)$ the vector $\psi^*$ is still unique. For instance, in the case in the case when the right side of an inclusion is represented in the form of intersection of two sets, the property required will be guaranteed (as explained in the example above) if $z$ is not a point equidistant from these sets. 
\end{rmrk}

\begin{prblm}
An interesting problem for future research is to explore the differential properties of the functional $\widetilde I(x,z)$ in the case of nonuniqueness of the vector $\psi^*$. The hypothesis is that the functional remains quasidifferentiable in this case as a superposition of the max-function and a quasidifferentiable function. Also note that in the case of intersection of two sets in the right-hand side of an inclusion the equality of the functions under min-function (at the same value of vector $\psi^*$) means that the points at which the minimum distances from the point $z$ to each of these sets is achieved, coincide, hence this is a distance from the point $z$ to their intersection which in turn means that the vector $\psi^*$ is unique in this case and the results of the functional $\widetilde I(x,z)$ subdifferentiability of this paper remain valid here.
\end{prblm}


Consider the model of double pendulum, one has
$$\ddot x_1 = -\nu^2(1+2 \mu^2) x_1 + \nu^2 \mu^2 x_2,$$
$$\ddot x_2 = \nu^2 x_1 - \nu^2 x_2, $$
where $\nu = \frac{m_1}{m_2}$, $\mu = \frac{g}{l}$ and $m_1$ and $m_2$ are the masses of the first and the second pendulum respectively, $g$ is the gravitation acceleration and $l$ is the thread length. In the literature (see, e. g., \cite{pend}) such models are considered with the assumption that the oscillations are small, therefore one can expect $\sin x_1$, $\sin x_2$ presented in the exact model to be approximately equal to $x_1$, $x_2$ respectively. Here it is suggested to consider a more precise model such that, for instance, the segments $[0.5 x_1, 1.5 x_1]$ and $[0.5 x_2, 1.5 x_2]$ are presented in the model instead including the $\sin x_1$ and $\sin x_2$ values (at least at the small values within the segment $[0, 0.5]$ of the \linebreak coordinates $x_1$, $x_2$). With such an assumption made one can note that the support function of the right-hand side is as follows:
$$c(F_1(x), \psi_1) = -\nu^2(1+2\mu^2) x_1 \psi_1 + 0.5 \nu^2(1+2\mu^2) |x_1| |\psi_1| + \nu^2 \mu^2 x_2 \psi_1 + 0.5 \nu^2 \mu^2 |x_2| |\psi_1|,$$
$$c(F_2(x), \psi_2) = \nu^2 x_1 \psi_2 + 0.5 \nu^2 |x_1| |\psi_2| - \nu^2 x_2 \psi_2 + 0.5 \nu^2 |x_2| |\psi_2|.$$

Consider the control system with disturbance:
$$
\dot x(t) = Ax(t) + B u(t) + \varpi (t)
$$

On the right side there is a control function $u(t)$, as well as a function $\varpi (t)$ playing the role of some disturbance in the system. Let, based on some physical considerations, the control as well as the perturbation $\varpi(t)$ belong to the sets depending on the phase coordinates  and suppose that these sets support functions have the structure discussed in the paper.  

Note that in book \cite{kurzh} an analytic solution to this problem is given using the apparatus of support functions as well as a conjugate function if one assumes that the corresponding concave hull of the sets required is realizable. This solution is obtained in the case when the control as well as the disturbance functions belong to some convex compact sets. Such analytic theory was possible due to the known formulas of the linear systems solution (via fundamental matrix) and support functions for standard convex compact sets of rather simple structure.

Consider the control system 
\begin{equation}
\label{cs}
\dot x(t) = Ax(t) + B u(t) 
\end{equation}
 on the time interval $[-t^*, T]$ (here $T$ is a given final time moment; see comments on the time moment~$t^*$ below). 

Let also ``discontinuity'' surface be given:
\begin{equation}
\label{ds}
s(x) = {\bf 0_m}
\end{equation}
where $s(x)$ is a known continuously differentiable $m$-dimensional vector-function. 

Consider the following form of controls:
\begin{equation}
\label{cv}
u_i(x) = -a_i |x| \mathrm{sgn}(s_i(x)), \ i = \overline{1, m}, \quad u_i = 0, \ i = \overline{m+1,n}
\end{equation} 
where $a_i \in [\underline {a}_i, \overline {a}_i]$, $i = \overline{1, m}$, are some positive numbers which are sometimes called gain factors.

In book \cite{Utkin} it is shown that if surface (\ref{ds}) is a hyperplane, then under natural assumptions and with sufficiently big values of the factors $a_i$, $i = \overline{1,m}$, controls (\ref{cv}) ensure system (\ref{cs}) hitting a small vicinity of this surface (\ref{ds}) from arbitrary initial state in the finite time $t^*$ and further staying in this neighborhood with the fulfillment of the condition $s_i(x(t)) \rightarrow 0$, $i = \overline{1,m}$, at $t \rightarrow \infty$. Now we are interested in the behavior of the system on the ``discontinuity'' surface (on the time interval $[0,T]$). 

We see that the right-hand sides of the first $m$ differential equations in system (\ref{cs}) with controls  (\ref{cv}) are discontinuous on the surfaces $s_i(x) = 0$, $i = \overline{1,m}$. So one has to use one of the known definitions of a discontinuous system solution. Thus, accordingly to Filippov definition \cite{fil} a solution of the discontinuous system considered satisfies differential inclusion
$$
\dot x_i \in A_i x + [\underline a_i, \overline a_i] |x| [-1, 1] = A_i x + [-\overline a_i, \overline a_i] |x| =: F_i(x), \ i = \overline {1,m}, \quad \dot x_i = A_i x, \ i = \overline{m+1,n},
$$
and we have
$$ c(F_i(x),\psi_i) = \psi_i A_i x + \overline a_i |x| |\psi_i| = \psi_i A_i x + \sum_{i=1}^n\max\{\overline a_i x_i \psi_i, -\overline a_i x_i \psi_i\}, \ i = \overline {1,m}, \quad c(F_i(x),\psi_i) = \psi_i A_i x, \ i = \overline{m+1,n}.$$

Consider the simpest physical model with dry friction. Let there is a block on a table subject to Coulomb friction on the contacting surface,
pulled by a force of gravity accelaration $g$. The differential equation \cite{stewart2009} for this system with the velocity $v$ is
\begin{equation} \label{coloumbfr} m \frac{dv(t)}{dt} = - \mu F_N \mathrm {Sgn}(v(t)) + g \end{equation}
where $\mathrm {Sgn}$ is a set-valued function given by
$$ \mathrm {Sgn}(v) = \left\{
\begin{array}{lll}
& \{1\}, \ v > 0, \\
& [-1, 1], \ v = 0, \\
& \{-1\}, \ v < 0.
\end{array}
\right.
$$

The quantity $F_N$ is the normal contact force (equal to $mg$ for a block of mass $m$) and $\mu$ is the
coefficient of Coulomb friction. We see that the right-hand sides in system (\ref{coloumbfr}) are discontinuous on the surfaces $s(v) = v = 0$. Here the approach is applied based on the Filippov works \cite{fil} on discontinuos differential equations theory, hence the system satisfies differential inclusion
$$
\dot v \in \mu F_N [-1, 1] + g,
$$

In a more general case a so called friction cone $FC$ is introduced and the equations of motion under the Coulomb dry friction law are of the form
$$M(q) \frac{dv}{dt} \in k(q, v) + FC(q)  $$ 
where $M(q)$ is a mass matrix and $k(q,v)$ is a function obtained from the kinetic and potential energy when applying Lagrange equations with generalized coordinates $q$ and velocities $v$. One should also impose holonomic or nonholonomic (unliteral) constraints as well as some complementary conditions considering: 1) the case when friction vector is on the boundary of the friction cone in the case of nonzero relative velocity, 2) physical situation where sliding stops during the contact period. The details may be found in \cite{Stewart}. The direct discretization is then applied there and the corresponding finite dimensional problem is solved via methods of linear algebra.

\begin{rmrk} Note that if from physical meaning the velocity itself is a discontinuous function of time, strictly speaking there is no sense in taking its derivative (so the acceleration is formally undefined). In order to overcome this difficulty, a possible generalization of a derivative as a measure and a relative concept of Radon-Nikodim derivative is used in some literature \cite{Moro} and the corresponding so called measure differential inclusions with such a derivative. The direct discretization is then applied there in order to get the corresponding finite dimensional problem. It is interesting to note: since the derivative there is a measure or a generalized function then in many cases it would be impossible to work with its discrete analogue; the solution suggested in \cite{Moro} is to consider the integral equation instead of the initial differential one and to ``discretize'' the function itself instead of its derivative in the left-hand side (and the integral of the function in the right-hand side) of the corresponding integral equation. This approach succeeds, for example, when we consider the Heaviside step-function and its generalized derivative, the delta-function. However, this paper does not deal with such cases.
\end{rmrk}

\section{The $\delta$-Quasidifferential Descent Method}

Describe the quasidifferentiable descent method for finding stationary points of the functional $I(x, z)$. 

Fix an arbitrary initial point $(x_{(1)}, z_{(1)}) \in C_n[0, T] \times P_n[0, T]$. Let the point $(x_{(k)}, z_{(k)}) \in C_n[0, T] \times P_n[0, T]$ be already constructed. If for each $t \in [0, T]$ minimum condition (\ref{11.400}) is satisfied (in practice, at discrete time moments $t_i$, $i = \overline{1,N}$, with some fixed accuracy $\overline{\varepsilon}$ in sense of Hausdorff norm in the space $R^n$, with some fixed discretization rank $N$), then $(x_{(k)}, z_{(k)})$ is a stationary point of the functional $I(x, z)$ and the process terminates. Otherwise, put
$$
(x_{(k+1)}, z_{(k+1)}) = (x_{(k)}, z_{(k)}) + \gamma_{(k)} G(x_{(k)}, z_{(k)})
$$
where the vector-function $G(x_{(k)}, z_{(k)})$ is the quasidifferential descent direction of the functional $I(x, z)$ at the point $(x_{(k)}, z_{(k)})$ and the value $\gamma_{(k)}$ is a solution of the following one-dimensional problem
\begin{equation}
\label{180}
\min_{\gamma \geq 0} I \Big( (x_{(k)}, z_{(k)}) + \gamma G(x_{(k)}, z_{(k)}) \Big) = I \Big( (x_{(k)}, z_{(k)}) + \gamma_{(k)} G(x_{(k)}, z_{(k)})\Big). 
\end{equation}
In practice, the problem above is solved on the interval $[0, \overline{\gamma}]$ with some  fixed $\overline \gamma$ value. Then one has the inequality $$I (x_{(k+1)}, z_{(k+1)}) < I (x_{(k)}, z_{(k)}) .$$ 

As seen from the algorithm described, in order to realize the $k$-th iteration, one has to solve four subproblems. The first subproblem is to calculate the quasidifferential of the functional $I(x, z)$ at the point $(x_{(k)}, z_{(k)})$. With the help of quasidifferential calculus rules the solution of this subproblem is obtained in formulas (\ref{4.44}), (\ref{4.6}), (\ref{4.55}), (\ref{11.2001}) and ~(\ref{11.2000}). \newline The second subproblem is to find the quasidifferential descent direction $G(x_{(k)}, z_{(k)})$; the following two paragraphs are devoted to solving this subproblem. \newline The third subproblem is one-dimensional minimization~(\ref{180}); there are many effective methods  (see, e. g.\cite{Vasil'ev}) for solving this subproblem. \newline The fourth subproblem is finding the values $\psi_i^*(x_{(k)}, z_{(k)})$, $i = \overline{1,n}$ at each time moment $t_i, i = \overline{1,N}$. Solving this problem is not straightforward since in general case it is of nonlinear optimization. However, there also many methods (see, e. g.\cite{Vasil'ev}) of nonlinear programming; moreover for some particular structures of the sets of the righ-hand sides of a differential inclusion it can be solved analytically.  

In order to obtain the vector-function $G(x_{(k)}, z_{(k)} )$, consider the problem
\begin{equation}
\label{19}
\max_{w \in \overline{\partial} I(x_{(k)}, z_{(k)}) } \min_{v \in \underline{\partial} I(x_{(k)}, z_{(k)})  } \int_0^T \big( v(t) + w(t) \big )^2 dt. 
\end{equation}
Denote $\overline{v}(t)$, $\overline{w}(t)$ its solution. (The vector-functions $\overline{v}(t)$, $\overline{w}(t)$, of course, depend on the point $(x_{(k)}, z_{(k)})$ but we omit this dependence in the notation for brevity.) Then the vector-function $G(x_{(k)}, z_{(k)}) = -\big( \overline{v} + \overline{w} \big)$ is a quasidifferential descent direction of the functional $I(x, z)$ at the point $(x_{(k)}, z_{(k)})$. 




In \cite{Fom_IJC} it is shown that 
 $$ \max_{w \in \overline{\partial} I(x_{(k)}, z_{(k)}) } \min_{v \in \underline{\partial} I(x_{(k)}, z_{(k)} ) } \int_0^T \big( v(t) + w(t) \big )^2 dt = $$
 \begin{equation} \label{eq}
  = \int_0^T \max_{w(t) \in \overline{\partial} \mathcal{I} (\xi_{(k)}(t), t) } \min_{v(t) \in \underline{\partial} \mathcal{I} (\xi_{(k)}(t), t) } \big( v(t) + w(t) \big )^2 dt. 
  \end{equation}
  
  The equality (\ref{eq}) justifies that in order to solve problem (\ref{19}) it is sufficient to solve the problem 
\begin{equation} \label{maxmin}
\max_{w(t) \in \overline{\partial} \mathcal{I} (\xi_{(k)}(t), t) } \min_{v(t) \in \underline{\partial} \mathcal{I} (\xi_{(k)}(t), t) } \big( v(t) + w(t) \big )^2
\end{equation}
for each time moment $t \in [0, T]$. Once again we emphasize that this statement holds true due to the special structure of the quasidifferential which in turn takes place due to the separation implemented of the vector functions \linebreak $x(t)$ and $\dot x(t)$ into ``independent'' variables.

Problem (\ref{maxmin}) at each fixed $t \in [0, T]$ is a finite-dimensional problem of finding the Hausdorff deviation of one convex compact set (a minus superdifferential) from another convex compact set (a subdifferential). This problem may be effectively solved in our case; its solution is described in the next paragraph. In practice, one makes a (uniform) partition of the interval $[0, T]$ and this problem is being solved for each point of the partition, i. e. one calculates $G((x_{(k)}, z_{(k)} ), t_i)$ where $t_i \in [0, T]$, $i = \overline{1, N}$, are discretization points (see notation of Lemma~6.1 below). Under additional natural assumption Lemma 6.1 below guarantees that the vector-function obtained via piecewise-linear interpolation of the quasidifferential descent directions calculated at each point of such a partition of the interval $[0, T]$ converges in the space $L^2_{2n}[0, T]$ (as the discretization rank~$N$ tends to infinity) to the vector-function $G(x_{(k)}, z_{(k)})$ sought. 

As noted in the previous paragraph, during the algorithm realization it is required to find the Hausdorff deviation of the minus superdifferential from the subdifferential of the functional $\mathcal {I} (x, z)$ at each time moment of a (uniform) partition of the interval $[0, T]$. In this paragraph we describe in detail a solution of this subproblem for some fixed value $t \in [0, T]$. It is known \cite{demvas} that in our case the subdifferential $\underline \partial \mathcal {I} (x, z, t)$ is a convex polyhedron $A(t) \subset R^{2n}$ and analogously the superdifferential $\overline \partial \mathcal {I} (x, z, t)$ is a convex polyhedron $B(t) \subset R^{2n}$. Herewith, of course, the sets $A(t)$ and $B(t)$ depend on the point $(x, z)$. For simplicity, we omit this dependence in this paragraph notation. Find the Hausdorff deviation of the set $-B(t)$ from the set $A(t)$. It is clear that in this case it is sufficient to go over all the vertices $b_j(t)$, $j = \overline{1,V}$ (here $V$ is a number of vertices of the polyhedron~$-B(t)$): find the Euclidean distance from every of these vertices to the polyhedron $A(t)$ and then among all the distances obtained choose the largest one. Let the Euclidean distance sought, corresponding to the vertex $b_j(t)$, $j = \overline{1,V}$, is achieved at the point $a_j(t) \in A(t)$ (which is the only one since $A(t)$ is a convex compact). Then the deviation sought is the value $||b_{\overline j}(t) - a_{\overline{j}}(t)||_{R^{2n}}$, $\overline j \in \{1, \dots, V\}$. (Herewith, this deviation may be achieved at several vertices of the polyhedron $-B(t)$; in this case $b_{\overline j}(t)$ denotes any of them.)
Note that the arising problem of finding the Euclidean distance from a point to a convex polyhedron can be effectively solved by various methods (see, e. g. \cite{wolfe}).

Give a lemma \cite{Fom_IJC} with a rather simple condition which, on the one hand, is quite natural for applications and, on the other hand, guarantees that the function $\mathcal L(t)$ obtained via piecewise-linear interpolation of the sought function $\mathcal G \in L^{\infty}_1[0, T]$ converges to this function in the space $L^2_1[0, T]$. 

\begin{lmm}
 Let the function $\mathcal G \in L^{\infty}_1[0, T]$ satisfy the following condition: for every $\overline{\delta} > 0$ the function $\mathcal G(t)$ is piecewise continuous on the set $[0, T]$ with the exception of only the finite number of the intervals $\big(\overline t_1(\overline{\delta}), \overline t_2(\overline{\delta})\big), \dots, \big(\overline t_{r}(\overline{\delta}), \overline t_{r+1}(\overline{\delta})\big)$ whose union length does not exceed the number $\overline{\delta}$. 
 
 Choose a (uniform) finite splitting $t_1 = 0, t_2, \dots, t_{N-1}, t_N = T$ of the interval $[0, T]$ and calculate the values $\mathcal G(t_i)$, $i = \overline{1, N}$, at these points. Let $\mathcal L(t)$ be the function obtained with the help of piecewise linear interpolation with the nodes $(t_i, \mathcal G(t_i))$, $i = \overline{1, N}$. Then for each $\overline\varepsilon > 0$ there exists such a number $\overline{N}(\overline\varepsilon)$ that for every $N > \overline{N}(\overline\varepsilon)$ one has $||\mathcal L - \mathcal G||^2_{L^2_1[0,T]} \leq \overline\varepsilon$.
 \end{lmm} 

Note that when implementing the algorithm using the rules of quasidifferential calculus, those functions that are active only with some error $\delta$ (or $\delta$-active) are considered. Let's introduce the concept of $\delta$-active function. Fix some value $\delta > 0$. We call the function $f_{\overline i}$, $\overline i \in \{1, \dots, M\}$, a $\delta$-active one at the point $x_0 \in R$ if $\varphi_1(x_0) = \max\limits_{i = \overline{1, M}} f_i(x) - f_{\overline i} \leq \delta$. We call the function $g_{\overline j}$, $\overline j \in \{1, \dots, K\}$, a $\delta$-active one at the point $x_0 \in R$ if $\varphi_2(x_0) = \min\limits_{j = \overline{1, K}} g_j(x) - g_{\overline j} \geq -\delta$. Taking this error into account justifies the method name, and also is crucial when proving its convergence in some sense in one special case (see Theorem \ref{th6.2} below).

It is clear that in most practical examples, the subproblems of the algorithm are also solved approximately. Herewith, the accuracy parameters there depend on specific methods for solving these problems and are also selected in advance.

Consider now a particular case when $c(F_i(x),\psi_i) = c_1(F_i(x),\psi_i)$, $i = \overline{1,n}$. Herewith, denote the functional minimized as $I_1(x,z)$. 

Put
\begin{equation}
\label{6.10}
\Psi_1(x, z) = \min_{||g||_{L_n^2[0, T] \times L_n^2[0,T]}=1}  \frac{\partial I_1(x, z)}{\partial g} = \frac{\partial I_1(x, z)}{\partial G} .
\end{equation}

Then necessary minimum condition \cite{demvar} of the functional $I_1(x, z)$ at the point $(x^*, z^*)$ is $\Psi_1(x^*, z^*) \geq 0$. 

Give now some results on the convergence of the quasidifferential descent method as applied to the functional $I_1(x,z)$. For this we have to make a few additional assumptions.

Fix the initial point $(x_{(1)}, z_{(1)})$. Suppose that the set $$lev_{I_1(x_{(1)}, z_{(1)})}I_1(\cdot, \cdot) = \Big\{ (x, z) \in C_n[0,T] \times P_n[0,T] \ \big| \ I_1(x, z) \leq I_1(x_{(1)}, z_{(1)})\Big\}$$ is bounded in $L_n^2[0, T] \times L_n^2[0,T]$-norm (due to the arbitrariness of the initial point, in fact, one must assume that the set $lev_{I_1(x_{(1)}, z_{(1)})}I_1(\cdot, \cdot)$ is bounded for every initial point taken).

Introduce now the set family $\mathcal{I}_1$. At first, define the functional $I_{1,p}$, $p = \overline{1, \left(\prod_{j=1}^r {m(j)}\right)^n}$ as follows. Its integrand is the same as the functional $I_1$ one, but the maximum function $\max\big\{f_{i, j_1}(x) \psi_i, \dots, f_{i, j_{m(j)}}(x) \psi_i\big\}$, $j = \overline{1,r}$, is substituted for each $i = \overline{1,n}$ by only one of the functions $f_{i, j_1} \psi_i, \dots, f_{i, j_{m(j)}} \psi_i$, $j \in \{1, \dots, r\}$. Let the family $\mathcal{I}_1$ consist of sums of the integrals over the intervals of the time interval $[0, T]$ splitting for all possible finite splittings. Herewith, the integrand of each summand in the sum taken is the same as some functional $I_{1,p}$ one, $p \in \left\{1, \dots, \left(\prod_{j=1}^r {m(j)}\right)^n \right\}$.

Let for every point constructed by the method described the following assumption be valid: the interval $[0, T]$ may be divided into a finite number of intervals, in every of which for each $i =\overline{1, n}$ either $h_i(x_{(k)},z_{(k)}) =0$, or one (several) of the functions $\left\langle - \psi_i^* \frac{\partial f_{i, j_q}(x_{(k)})}{\partial x}, G_1(x_{(k)},z_{(k)}) \right \rangle$, $j = \overline{1,r}$, $q = \overline{1,m(j)}$, is (are) active.
\newpage

For functionals from the family $\mathcal I_1$ we make the following additional assumption. Let there exist such a finite number $L$ that for every $\hat I_1 \in \mathcal I$ and for all $\overline x, \overline z, \overline{\overline x}, \overline{\overline z}$ from a ball with the center in the origin and with some finite radius $\hat r + \hat \alpha$ (here $\hat r > \sup\limits_{(x,z) \in lev_{I_1(x_{(1)}, z_{(1)})}I_1(\cdot, \cdot)}||(x, z)||_{L^2_n[0,T] \times L^2_n[0,T] }$ and $\hat \alpha$ is some positive number) one has 
\begin{equation}
\label{6.100} ||\nabla \hat I_1(\overline x, \overline z) -  \nabla \hat I_1(\overline{\overline x}, \overline{\overline z})||_{L^2_n[0,T] \times L^2_n[0,T] } \leq L ||({\overline x}, {\overline z})' - (\overline{\overline x}, \overline{\overline z})' ||_{L^2_n[0,T] \times L^2_n[0,T]}. 
\end{equation}

\begin{thrm} \label{th6.2} 
Under the assumptions made one has the inequality
\begin{equation} 
\label{7_new_4}
\underline \lim_{k \rightarrow \infty} \Psi_2(x_{(k)}, z_{(k)}) \geq 0
\end{equation}
for the sequence built according to the rule above. (See \cite{Fom_cocv}.)
\end{thrm}

In \cite{Fom_cocv} the explanations of the assumptions above with illustrating examples are made as well.

Consider now a particular case when $c(F_i(x),\psi_i) = c_2(F_i(x),\psi_i)$, $i = \overline{1,n}$. Herewith, denote the functional minimized as $I_2(x,z)$. 

 
In order to simplify the presentation, we preliminarily explore the convergence properties of the method proposed for an analogous problem with the functional in ``standard'' form; and then turn to a more general problem considered in this paper. 

So consider the problem of minimization of a functional which has a maximum of the finite number of continuously differentiable functions as its integrand. Let the functional on the space $\Big(C_n[0,T] \times P_n[0,T], \linebreak || \cdot ||_{L^2_n [0, T]} \times || \cdot ||_{L^2_n [0, T]} \Big)$ be 
$$I_2(x,z) = \int_0^T \xi(x(t), z(t)) dt = \int_0^T \max\Big\{ \xi_1(x(t), z(t)), \dots, \xi_M(x(t),z(t)) \Big\} dt,$$
where $\xi_i(x, z)$, $i = \overline{1, M}$, are continuously differentiable on $\mathbb R^n \times \mathbb R^n$ functions. 

We suppose that the points $x_k(t)$ and $z_k(t)$ and the subdifferential descent directions $G(x_k, z_k, t)$ are bounded (in the uniform norm) uniformly in $k = 1, 2, \dots$ during the described method implementation.

Note that since the variables $x$ and $z$ are considered as independent ones here, in fact, the problem above may be solved by pointwise minimization of the integrand. However, we ignore this simpification and give a ``full-flegged'' proof in order to convey the technique, applicable to the more general case of the functional considered in this chapter.

Put
$$
\Psi_2(x, z) = \min_{||g||_{L_n^2[0, T] \times L_n^2[0,T]}=1}  \frac{\partial I_2(x, z)}{\partial g} = \frac{\partial I_2(x, z)}{\partial G} = 
$$
$$ = \int_0^T \max_{i \in \mathcal R(x(t), z(t))} \left( \left \langle \frac{\partial \xi_i(x(t), z(t))}{\partial x}, G_1(x, z, t) \right \rangle + \left \langle \frac{\partial \xi_i(x(t), z(t))}{\partial z}, G_2(x, z, t) \right \rangle \right) dt,$$
here $\mathcal R(x, z) = \Big\{i \in \{1, \dots, M\} \ | \ \xi_i(x, z) = \xi(x, z) \Big\}$ and $G_1(x,z) (G_2(x,z))$ consists of the first (last) $n$ components of the vector-function $G(x, z)$.

Then necessary minimum condition \cite{demvar} of the functional $I_2(x, z)$ at the point $(x^*, z^*)$ is $\Psi_2(x^*, z^*) \geq 0$. 

Let
\begin{equation}
\label{11.55}
\Psi_{2,\overline\Delta}(x, z) = \int_0^T \max_{i \in \mathcal R_{\overline\Delta}(x(t), z(t))} \left( \left \langle \frac{\partial \xi_i(x(t), z(t))}{\partial x}, G_1(x, z, t) \right \rangle + \left \langle \frac{\partial \xi_i(x(t), z(t))}{\partial z}, G_2(x, z, t) \right \rangle \right) dt,
\end{equation}
here $\mathcal R_{\overline\Delta}(x, z) = \Big\{i \in \{1, \dots, M\} \ | \ \xi(x, z) - \xi_i(x, z) \leq {\overline\Delta} \Big\}$.

If one has $\Psi_{2,\overline\Delta}(x^{\overline\Delta *}, z^{\overline\Delta *}) \geq 0$, then the point $(x^{\overline\Delta *}, z^{\overline\Delta *})$ is called a $\overline\Delta$-stationary point of the functional $I_2(x, z)$.

The concept of a $\overline\Delta$-stationary point is as follows. Arguing as in \cite{demmal} one can show that if the functions $\xi_i(x,z)$ are convex at $i \in \mathcal R_{\overline\Delta}(x^{\overline\Delta *}, z^{\overline\Delta *})$, then one has the inequality $$0 \leq I_2(x^{\overline\Delta *}, z^{\overline\Delta *}) - \min_{(x, z) \in  C_n[0,T] \times P_n[0,T]} \int_0^T \xi(x(t), z(t)) dt  \leq T \overline\Delta,$$ i. e. the functional $I_2$ value at its $\overline\Delta$-stationary point does not differ much from the functional $I_2$ value at its minimum point.  

\begin{wrapfigure}{r}{0.35\textwidth} \begin{center}
     \includegraphics[width=0.35\textwidth]{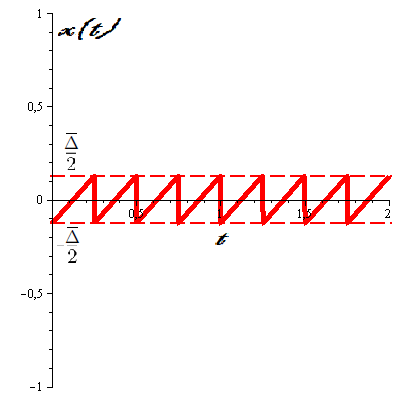}
  \end{center}    
      \end{wrapfigure}
Note that considering the case of the max-function $ \xi(x(t), z(t))$, $t \in [0, T]$, under the integral we make a modification into the algorithm described. The modification is as follows: when calculating the subdifferential of this max-function at each time moment we consider the function $\xi_{\overline i}(x, z)$, $\overline i \in \{1,\dots,M\}$, to be active at the point $(x_0, z_0) \in \mathbb R^n \times \mathbb R^n$ if it is active with an error of $\overline\Delta$. i. e. $\ \xi(x_0, z_0) - \xi_{\overline i}(x_0, z_0) \leq \overline\Delta$. In other words, when using the corresponding subdifferential calculus rule during the method implementation we use the set $R_{\overline\Delta}(x_k(t_i), z_k(t_i))$ instead of the set $R(x_k(t_i), z_k(t_i))$ at each time moment $t_i$, $i = \overline{1,N}$, of the time interval $[0, T]$ discretization.
\bigskip

Let us illustrate this modification by a simple example. Consider the minimization of the functional \bigskip

$$ J(x) = \int_0^2 \xi{(x(t))} dt = \int_0^2 \max\{\xi_1(x(t)), \xi_2(x(t))\} dt = \int_0^2 \max\{x(t), -x(t)\} dt = \int_0^2 |x(t)| dt.$$
\bigskip
\bigskip

Obviously $x^*(t) = 0$ for all $t \in [0,T]$. Take the point $x_1(t) = t-1$ as an initial approximation. During the modified method implementation, at some step the ``saw'' curve showed in the picture will be obtained, herewith it will not exceed the prescibed value $\overline\Delta$ at each $t \in [0,2]$. So if one considers the $\overline\Delta$-active functions at any discretization point $t_i$, $i = \overline{1,N}$, then the minimum condition is satisfied since $0 \in [-1, 1] =$ $= co\Big\{\frac{\partial \xi_1(x^{\overline\Delta *}(t_i))}{\partial x}, \frac{\partial \xi_2(x^{\overline\Delta *}(t_i))}{\partial x}\Big\}$, and the process terminates at this step. The $\overline\Delta$-stationary point, the ``saw'' curve $x^*(t)$ obtained, is an approximate solution and $0 \leq J(x^{\overline\Delta *}) - J(x^*) \leq T \overline\Delta = 2 \overline\Delta$.

\begin{thrm} \label{th11.5.3} 
Under the assumptions made one has the inequality
$$
\underline \lim_{k \rightarrow \infty} \Psi_{2, \overline \Delta} (x_{(k)}, z_{(k)}) \geq 0
$$
for the sequence built according to the (modified) rule above. 
\end{thrm}
\begin{proof}
Denote $$\psi_{2, \overline \Delta}(x(t), z(t)) = \max_{i \in R_{\overline\Delta}(x(t), z(t))} \left( \left \langle \frac{\partial \xi_i(x(t), z(t))}{\partial x}, G_1(x, z, t) \right \rangle + \left \langle \frac{\partial \xi_i(x(t), z(t))}{\partial z}, G_2(x, z, t) \right \rangle \right), $$
$$C_1 = \max_{t \in [0,T]}\max_{(x_{k}(t), z_{k}(t)) \in \mathbb R^n \times \mathbb R^n} \max_{i = \overline{1,M}} \Big|\Big| \Big( \frac{\partial \xi_i(x_{k}(t), z_{k}(t))}{\partial x}, \frac{\partial \xi_i(x_{k}(t), z_{k}(t))}{\partial z} \Big)' \Big|\Big|_{\mathbb R^n \times \mathbb R^n},$$
$$C_2 = \max_{t \in [0,T]}\max_{(x_{k}(t), z_{k}(t)) \in \mathbb R^n \times \mathbb R^n} \psi_{2, \overline \Delta}(x_k(t), z_k(t))$$
 and note that $C_1 < \infty$ and $C_2 < \infty$ by assumption.

Assume the contrary. Then there exist such a subsequence $\{(x_{k_j},z_{k_j})\}_{j=1}^{\infty}$ and such a number $b > 0$ that for each $j  \in \mathbb N$ we have the inequality
$$
\int_T \psi_{2, \overline \Delta}(x_{k_j}(t), z_{k_j}(t)) dt \leq -b.
$$

Then due to the boundedness of the integrand $\psi_{2, \overline \Delta}(x(t), z(t))$ on $[0,T]$ it is clear that there exists such a number $\overline b > 0$ that for each $j  \in \mathbb N$ one has the inequality
\begin{equation}
\label{11.55}
 \int_{T^-_{k_j}} \psi_{2, \overline \Delta}(x_{k_j}(t), z_{k_j}(t)) dt + \int_{T^+_{k_j}} \psi_{2, \overline \Delta}(x_{k_j}(t), z_{k_j}(t)) dt  \leq -\frac{3 b}{4},
 \end{equation}
where the subsets $T^-_{k_j}$ and $T^+_{k_j}$ of $[0,T]$ are defined as follows:
$$T^-_{k_j} = \Big\{ t \in [0,T] \ \big| \ \psi_{2, \overline \Delta}(x_{k_j}(t), z_{k_j}(t)) < -\overline b \Big\},$$
 $$T^+_{k_j} = \Big\{ t \in [0,T] \ \big| \ \psi_{2, \overline \Delta}(x_{k_j}(t), z_{k_j}(t)) > 0 \Big\},$$
 also introduce the set
 $$T^0_{k_j} = \Big\{ t \in [0,T] \ \big| \ -\overline b \leq \psi_{2, \overline \Delta}(x_{k_j}(t), z_{k_j}(t)) \leq 0 \Big\}.$$
 
Below we use the facts that the functions $\xi_i(x, z)$, $i = \overline{1,M}$, are continuously differentiable and the sequences $\{(x_{k_j}(t), z_{k_j}(t))\}_{j = 1}^{\infty}$ and $\{G_{k_j}(t)\}_{j = 1}^{\infty}$, $G_{k_j}(t) := G(x_{k_j}, z_{k_j}, t)$, $j \in \mathbb N$, are bounded uniformly on $[0, T]$ by assumption.
 
 Consider the case $t \in T^0_{k_j}$. 
 \newline  
 For each $j \in \mathbb N$, if $i \in R_{\overline\Delta}(x_{k_j}(t), z_{k_j}(t))$, then for each $\alpha > 0$ one has 
 $$\xi_i\Big((x_{k_j}(t), z_{k_j}(t)) + \alpha G_{k_j}(t)\Big) \leq \xi_i(x_{k_j}(t), z_{k_j}(t)) + \alpha \cdot 0 + o_i \Big(\alpha, x_{k_j}(t), z_{k_j}(t), G_{k_j}(t)\Big).$$ 
There exists \cite{demmal} the number $\overline{\overline{\overline \alpha}}$, which does not depend on the number $k_j$ and on the index $i$, such that for each $\alpha \in (0, \overline{\overline{\overline \alpha}}]$ one has the inequality
 $$\xi_i\Big((x_{k_j}(t), z_{k_j}(t)) + \alpha G_{k_j}(t)\Big) \leq \xi_i(x_{k_j}(t), z_{k_j}(t)) + {{ \alpha}} \frac{b}{4 T}.$$
 For each $j \in \mathbb N$, if $i \notin R_{\overline\Delta}(x_{k_j}(t), z_{k_j}(t))$, then for each $\alpha > 0$ one has the inequality
 $$\xi_i\Big((x_{k_j}(t), z_{k_j}(t)) + \alpha G_{k_j}(t)\Big) \leq \xi_i(x_{k_j}(t), z_{k_j}(t)) + \alpha C_1 + o_i \Big(\alpha, x_{k_j}(t), z_{k_j}(t), G_{k_j}(t)\Big) \leq$$
 $$ \leq \xi(x_{k_j}(t), z_{k_j}(t)) - \overline \Delta + \alpha C_1 + o_i \Big(\alpha, x_{k_j}(t), z_{k_j}(t), G_{k_j}(t)\Big).$$
Since $C_1 < \infty$, there exists \cite{demmal} the number $\overline{\overline{\overline \alpha}}_0$, $0 < \overline{\overline{\overline \alpha}}_0 \leq \overline{\overline{\overline \alpha}}$, which does not depend on the number $k_j$ and on the index $i$, such that for each $\alpha \in (0, \overline{\overline{\overline \alpha}}_0]$ one has
 $$
 \xi_i\Big((x_{k_j}(t), z_{k_j}(t)) + \alpha G_{k_j}(t)\Big) \leq \xi_i(x_{k_j}(t), z_{k_j}(t)) - \frac{\overline \Delta}{2}.
$$
 
 For each $\alpha \in (0, \overline{\overline{\overline \alpha}}_0]$, since $\max \Big \{ {{ \alpha}} \frac{b}{4T}, -\frac{\overline \Delta}{2} \Big\} = {{ \alpha}} \frac{b}{4T}$, we henceforth obtain the inequality
 \begin{equation}
 \label{11.56}
  \int_{T^0_{k_j}}\xi\Big((x_{k_j}(t), z_{k_j}(t)) + \alpha G_{k_j}(t)\Big)  dt \leq \int_{T^0_{k_j}}\xi(x_{k_j}(t), z_{k_j}(t)) dt + {{ \alpha}} \frac{b}{4}.
   \end{equation}
 
 Consider the case $t \in T^-_{k_j}$.
 \newline
 If $i \notin R_{\overline\Delta}(x_{k_j}(t), z_{k_j}(t))$, then (since $C_1 < \infty$) there exists \cite{demmal} the number ${\overline{\overline \alpha}}_0$, $0 < {\overline{\overline \alpha}}_0 \leq \overline{\overline{\overline \alpha}}_0$, which does not depend on the number $k_j$ and on the index $i$, such that for each $\alpha \in (0, {\overline{\overline \alpha}}_0]$ one has 
 $$\xi_i\Big((x_{k_j}(t), z_{k_j}(t)) + \alpha G_{k_j}(t)\Big) \leq \xi_i(x_{k_j}(t), z_{k_j}(t)) - \frac{\overline \Delta}{2}.$$
   For each $j \in \mathbb N$, if $i \in R_{\overline\Delta}(x_{k_j}(t), z_{k_j}(t))$, then for each $\alpha > 0$ one has  
 $$
 \xi_i\Big((x_{k_j}(t), z_{k_j}(t) + \alpha G_{k_j}(t)\Big)) \leq \xi_i(x_{k_j}(t), z_{k_j}(t)) + \alpha \psi_{2, \overline \Delta}(x_{k_j}(t), z_{k_j}(t)) + o_i \Big(\alpha, x_{k_j}(t), z_{k_j}(t), G_{k_j}(t)\Big).
 $$
Since $C_2 < \infty$, there exists \cite{demmal} the number ${\overline{\overline \alpha}}$, $0 < {\overline{\overline \alpha}} \leq {\overline{\overline \alpha}_0}$, which does not depend on the number $k_j$ and on the index $i$, such that for each $\alpha \in (0, {\overline{\overline \alpha}}]$ one has the inequality
 \begin{equation}
 \label{11.57}
 \xi_i\Big((x_{k_j}(t), z_{k_j}(t)) + \alpha G_{k_j}(t)\Big) \leq \xi_i(x_{k_j}(t), z_{k_j}(t)) + {{ \alpha}} \psi_{2, \overline \Delta}(x_{k_j}(t), z_{k_j}(t)) + {{ \alpha}} \frac{b}{8 T}
 \end{equation}
 and besides due to the inequality $\overline \Delta > 0$ and the boundedness of the function $\psi_{2, \overline \Delta}(x_{k_j}(t), z_{k_j}(t))$ on $[0,T]$ and due to the set $T^-_{k_j}$ definition, this number ${\overline{\overline \alpha}}$ may be taken such that for each $\alpha \in (0, {\overline{\overline \alpha}}]$ it will be $\max \Big\{ {{ \alpha}} \psi_{2, \overline \Delta}(x_{k_j}(t), z_{k_j}(t)), -\frac{\overline \Delta}{2} \Big\} = {{ \alpha}} \psi_{2, \overline \Delta}(x_{k_j}(t), z_{k_j}(t)).$
 
 Consider the case $t \in T^+_{k_j}$.
 \newline
 If $i \notin R_{\overline\Delta}(x_{k_j}(t), z_{k_j}(t))$, then (since $C_1 < \infty$) there exists \cite{demmal} the number ${{\overline \alpha}}_0$, $0 < {{\overline \alpha}}_0 \leq {\overline{\overline \alpha}}$, \linebreak which does not depend on the number $k_j$ and on the index $i$, such that for each $\alpha \in (0, {{\overline \alpha}}_0]$ one has
 $$\xi_i\Big((x_{k_j}(t), z_{k_j}(t)) + \alpha G_{k_j}(t)\Big) \leq \xi_i(x_{k_j}(t), z_{k_j}(t)) - \frac{\overline \Delta}{2}.$$
  For each $j \in \mathbb N$, if $i \in R_{\overline\Delta}(x_{k_j}(t), z_{k_j}(t))$, then for each $\alpha > 0$ one has 
 $$\xi_i\Big((x_{k_j}(t), z_{k_j}(t)) + \alpha G_{k_j}(t)\Big) \leq \xi_i(x_{k_j}(t), z_{k_j}(t)) + \alpha \psi_{2, \overline \Delta}(x_{k_j}(t), z_{k_j}(t)) + o_i \Big(\alpha, x_{k_j}(t), z_{k_j}(t), G_{k_j}(t)\Big).$$ 
 Since $C_2 < \infty$, there exists \cite{demmal} the number ${{\overline \alpha}}$, $0 < {{\overline \alpha}} \leq {{\overline \alpha}_0}$, which does not depend on the number $k_j$ and on the index $i$, such that for each $\alpha \in (0, {{\overline \alpha}}]$ one has the inequality
 \begin{equation}
 \label{11.58} \xi_i\Big((x_{k_j}(t), z_{k_j}(t)) + \alpha G_{k_j}(t)\Big) \leq \xi_i(x_{k_j}(t), z_{k_j}(t)) + {{\alpha}} \psi_{2, \overline \Delta}(x_{k_j}(t), z_{k_j}(t)) + {{\alpha}} \frac{b}{8 T}
 \end{equation}
 and note that due to the set $T^+_{k_j}$ definition for each $\alpha \in (0, {{\overline \alpha}}]$ it will be $\max \Big\{ {{\alpha}} \psi_{2, \overline \Delta}(x_{k_j}(t), z_{k_j}(t)), -\frac{\overline \Delta}{2} \Big\} =$ $={{\alpha}} \psi_{2, \overline \Delta}(x_{k_j}(t), z_{k_j}(t)).$
 \bigskip
 
 From (\ref{11.55}), (\ref{11.56}), (\ref{11.57}), (\ref{11.58}) we finally have 
 $$ \int_{T}\xi\Big((x_{k_j}(t), z_{k_j}(t)) + \overline \alpha G_{k_j}(t)\Big)  dt = $$ 
 {\small $$ = \int_{T^-_{k_j}}\xi\Big((x_{k_j}(t), z_{k_j}(t)) + \overline \alpha G_{k_j}(t)\Big) dt + \int_{T^+_{k_j}}\xi\Big((x_{k_j}(t), z_{k_j}(t)) + \overline \alpha G_{k_j}(t)\Big)  dt +$$ $$+ \int_{T^0_{k_j}}\xi\Big((x_{k_j}(t), z_{k_j}(t)) + \overline \alpha G_{k_j}(t)\Big) dt \leq $$ }
 $$\leq \int_{T^-_{k_j}}\xi(x_{k_j}(t), z_{k_j}(t)) dt + \int_{T^-_{k_j}} \overline \alpha \psi_{2, \overline \Delta}(x_{k_j}(t), z_{k_j}(t)) dt + \overline \alpha \frac{b}{8} + $$
  $$+ \int_{T^+_{k_j}}\xi(x_{k_j}(t), z_{k_j}(t)) dt + \int_{T^+_{k_j}} \overline \alpha \psi_{2, \overline \Delta}(x_{k_j}(t), z_{k_j}(t)) dt + \overline \alpha \frac{b}{8} + $$
 $$+ \int_{T^0_{k_j}}\xi(x_{k_j}(t), z_{k_j}(t)) dt + {{\overline \alpha}} \frac{b}{4} \leq $$
 $$ \leq \int_{T}\xi(x_{k_j}(t), z_{k_j}(t)) dt - \overline{\alpha} \frac{3 b}{4} + \overline{\alpha} \frac{b}{4} + \overline{\alpha} \frac{b}{4} =$$ $$= \int_{T}\xi(x_{k_j}(t), z_{k_j}(t)) dt - \overline{\alpha} \frac{b}{4} = \int_{T}\xi(x_{k_j}(t), z_{k_j}(t)) dt - \beta $$
 uniformly in $j \in \mathbb N$.

As one may directly check (using the set $R_{\overline\Delta}(x_{(k)},z_{(k)})$ definition), the sequence $\left\{I_2(x_{(k)},z_{(k)})\right\}_{k=1}^{\infty}$ is monotonically decreasing and bounded below by zero (recall that the functional $I(x,z)$ is nonnegative by construction), hence, it has a limit:
\label{11.new10000} 
$$
\left\{I_2(x_{(k)},z_{(k)})\right\} \rightarrow I_2^* \ \mathrm{at} \ k \rightarrow \infty,
$$
  herewith, at each $k \in \mathbb N$ one has $\left\{I_2(x_{(k)},z_{(k)})\right\} \geq I_2^*$. Further, the contradiction can be obtained similarly as it is done in the end of Theorem \ref{th6.2} proof. 
  \end{proof}
  
  Now turn back to the problem of the paper. Here we consider the case when $c(F_i(x),\psi_i) =$ $c_2(F_i(x),\psi_i)$, $i = \overline{1,n}$. Herewith, denote the functional minimized as $I_2(x,z)$. 
  
  Put
$$
\Psi_2(x, z) = \min_{||g||_{L_n^2[0, T] \times L_n^2[0,T]}=1}  \frac{\partial I_2(x, z)}{\partial g} = \frac{\partial I_2(x, z)}{\partial G}.
$$

Then necessary minimum condition \cite{demvar} of the functional $I_2(x, z)$ at the point $(x^*, z^*)$ is $\Psi_2(x^*, z^*) \geq 0$. 

The functional $\Psi_{2,\overline\Delta}(x, z)$ is constructed for the functional $I_2(x, z)$ by the formula analogous to (\ref{11.55}). 


Give now some results on the convergence of the quasidifferential descent method as applied to the functional $I_2(x,z)$. For this we have to make a few additional assumptions.
  
Let for every point constructed by the method described the following assumption be valid: the interval $[0, T]$ may be divided into a finite number of intervals, in every of which for each $i =\overline{1, n}$ either $h_i(x_{(k)},z_{(k)}) =0$, or one (several) of the functions $\left\langle - \psi_i^* \frac{\partial g_{i, j_p}(x_{(k)})}{\partial x}, G_1(x_{(k)},z_{(k)}) \right \rangle$, $j = \overline{1,s}$, $p = \overline{1,k(j)}$, is (are) active (where $G_1(x_{(k)},z_{(k)})$ consists of the first $n$ components of the vector-function $G(x_{(k)},z_{(k)})$).
  
  \begin{thrm} \label{th11.5.4} 
Under the assumptions made one has the inequality
$$
\underline \lim_{k \rightarrow \infty} \Psi_{2, \overline \Delta}(x_{(k)}, z_{(k)}) \geq 0
$$
for the sequence built according to the (modified) rule above. 
\end{thrm}

\begin{rmrk}
As is seen, although related from the optimization theory viewpoint, still slightly different assumptions are made regarding the functionals $I_1$ and $I_2$ properties and different techniques for the corresponding proofs are used. It is an interesting question to ``relate'' the assumptions and the proofs for both the results, since it seems that they may be ``unified''.
\end{rmrk}

\section{Numerical Examples}
Let us return to some examples noted in Section 5. In all the examples considered the values $\delta = 10^{-3}$, $\varepsilon = 10^{-2}$ and $\overline \gamma = 1$ of the parameters were taken in the $\delta$-quasidifferential descent method.

\begin{xmpl} Consider the differential inclusion
$$\dot x_1 \in [x_1, x_2] \cap [x_2-1, x_2 + 1],$$ 
$$\dot x_2 = x_1 + 1$$
on the time interval $[0, 1]$ with the boundary conditions
$$x_1(0) = 0, \, x_2(0) = 0, \quad x_2(1) = 1.$$
Consider the more general problem of minimizing the functional
$$\mathcal J(x) = \int_0^1 |x_1(t)| dt.$$

One of the obvious solutions is 
$x^*_1(t) = 0$, $x_2^*(t) = t$ for all $t \in [0,1]$.

Such an illustrative example is taken intentionally in order for the cost functional to be substantially nonsmooth (but only subdifferentiable) on the optimal trajectory. In order to solve this optimization problem with restrictions consider the correpsonding unconstrained problem of minimizing the functional 
$$
J(x, z) := \mathcal J(x) + \lambda I(x, z)
$$
with sufficiently big penalty factor $\lambda$ value.
 Also note that $ J(x^*, z^*) = 0$. 


We have
$$ c(F_1(x), \psi_1) = \min \left\{\frac{x_1+x_2}{2}\psi_1 + \frac{|x_1-x_2|}{2}, \, x_2 \psi_1 + 1\right\}, \quad c(F_2(x), \psi_2) = \psi_2(x_1+1)$$
for $\psi_1, \psi_2 \in S_1$.

Consider the support function $c(F_1(x), \psi_1)$ for $\psi_1 \in S_1$ in details. As proved in Section 4, when both functions are active, this function is superdifferentiable provided that $x_1 \neq x_2$. But it is interesting to note that in the case $x_1 = x_2$ the first function is the only one to be active, and since abs-function can be represented as the max-function of two continuous differentiable functions (the function under abs itself and the opposite one), then as again proved in Section 4 this function is subdifferentiable in the case $x_1 = x_2$.

Take $\lambda^* = 10$ and $(x_{(1)}, z_{(1)}) = (-1, -2, 0, 0)'$ as the first approximation. We intentionally take such a point in order to make both functions under minimum in the first support function to be active (see the formula above and note that, as one can check, $\psi_1^*(x_{(1)}, z_{(1)}) = 1$), hence the support function would be substantionally nonsmooth at this point. At the end of the process the discretization step was equal to $10^{-1}$. Figure~1 illustrates the trajectories obtained. From Figure 1 we see that the differential inclusion is practically satisfied (we see that the resulting curve $\dot x_1(t)$ practically coincides with the ``lower'' boundary $x_1(t)$ of the set $F_1(x(t))$ obtained and the resulting curve~$\dot x_2(t)$ practically coincides with the curve $x_1(t)+1$ (which is formally the set $F_2(x(t))$ obtained)). The boundary values error doesn't exceed the magnitude $5 \times 10^{-3}$. To obtain such an accuracy $37$ iterations have been required. The functional value on the trajectory obtained is of order $10^{-3}$.
\begin{figure*}[h!]
\begin{minipage}[h]{0.3\linewidth}
\center{\includegraphics[width=1\linewidth]{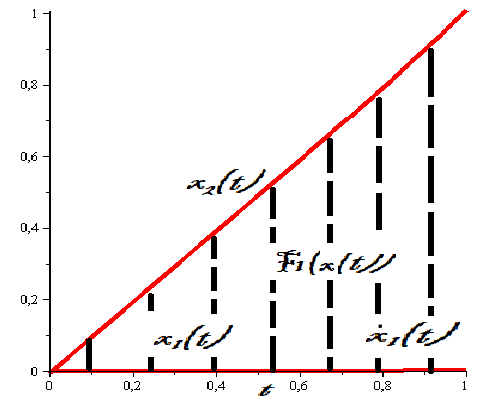} }
\end{minipage}
\hfill
\begin{minipage}[h]{0.3\linewidth}
\center{\includegraphics[width=1\linewidth]{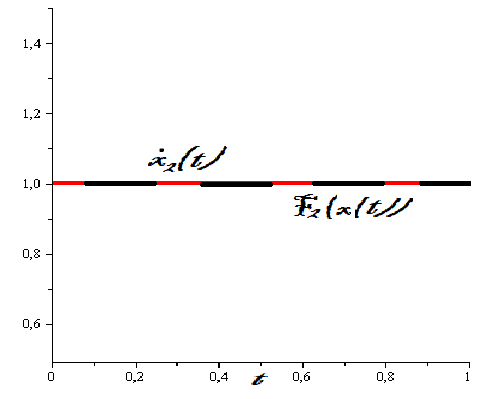} }
\end{minipage}
\caption{Solution of Example 7.1}
\end{figure*} 
\end{xmpl}

\begin{xmpl}
Consider the differential inclusion
$$
\dot{x} \in F(x), \quad F(x) =
B_2 \cap E_2(x), \quad t \in [0, 1],
$$
where $E_2(x)$ is an ellipse depending on the phase coordinate
$x_1$:
$$\displaystyle{E_2(x) = \Big\{ (\overline x, \overline y) \in R^2 \
\big| \ \frac{(\overline x - 2)^2}{x_1^2+3} + \overline{y}^2 \leq 1 \Big\} }.$$
The boundary conditions are 
$$
x(0) = (0, 0)', \quad (0.75, 0.5)'.
$$

We have
$$c(B_2, \psi) = \sqrt{\psi_1^2 + \psi_2^2} = 1, \quad c(E_2(x), \psi) = \sqrt{(x_1^2+3) \psi_1^2 + \psi_2^2} + 2 \psi_1.$$

Note that here we use the functional $\widetilde I(x,z)$ (see Remark \ref{rm5.1}).

Take an initial point $x_{(1)} = (t, t)'$ and $z_{(1)} = (1, 1)'$. At the end of the process the discretization step was equal to $10^{-1}$. Figure \ref{pic9.12} illustrated the trajectories obtained. In Figure \ref{pic9.12} the points $z^*(t)$ and the allowed set $B_2 \cap {E}_2(x(t))$ of these points location are depicted at some $t$-values from the segment $[0,1]$, it is seen that the differential inclusion considered is satisfied for these values $t \in [0,1]$ (it is easy to chek that it is correct for all the others time moments as well $t \in [0,1]$). The boundary values error doesn't exceed the magnitude $5 \times 10^{-3}$. To obtain such an accuracy $15$ iterations have been required. The functional $I(x^*, z^*)$ value on the trajectory obtained is of order $5 \times 10^{-3}$.

Note that as can be easily checked the distances from this point to the sets considered is not equal for each time moment $t \in [0,1]$, hence the requirement of the vector $\psi^*$ uniqueness (in order for the Theorem \ref{th4.2} to be valid) is fullfilled here (see Remark \ref{rm5.1}). The same is true for the point $z_{(k)}$, $k = 1,2, \dots$ on all the other iterations of the algorithm implemented.


\begin{figure*}[h!] \label{pic9.12}
\begin{minipage}[h]{0.26\linewidth}
\center{\includegraphics[width=1\linewidth]{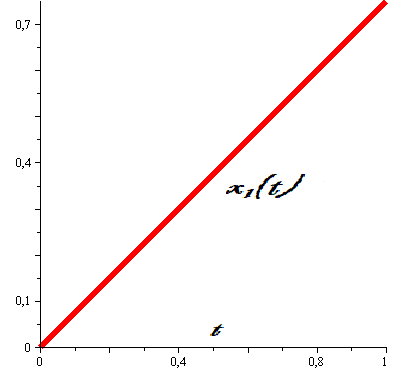} }
\end{minipage}
\hfill
\begin{minipage}[h]{0.26\linewidth}
\center{\includegraphics[width=1\linewidth]{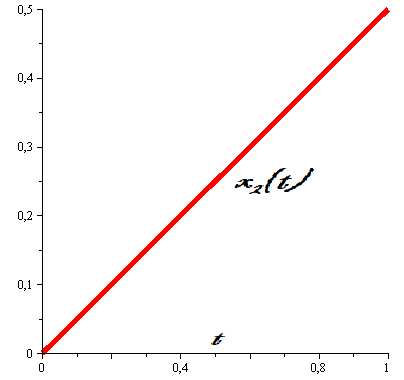} }
\end{minipage}
\begin{minipage}[h]{0.37\linewidth}
\center{\includegraphics[width=1\linewidth]{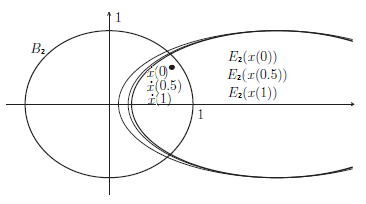} }
\end{minipage}
\caption{Solution of Example 7.2}
\label{ris:image2}
\end{figure*}
\end{xmpl}
\smallskip

\begin{xmpl}
The test problem used here is a modified problem from \cite{mult} with a system of three masses connected by springs with a forcing term and friction as shown in Figure 3. This object combines the problems of two types: modeling dry friction generates a discontinuous system moving in a sliding mode, thus can be considered as a differential inclusion; the sliding mode on the surface $s(x) = x_4 = 0$ is provided by an additional nonsmooth control function $u_1(x) = -a_1 |x| \mathrm{sgn}(s(x))$ where the gain factor $a_1 = 1.05$,  giving a support function of the set in the right-hand side of this inclusion in the form discussed in the paper. The aim of external forcing term $u(t)$ is to bring the system (moving in a sliding mode) to the point prescribed at the final moment of time.

The actual differential
equations to be solved are
$$\ddot x_1 = ( - x_1) + (x_2 - x_1) - \dot x_1 - 1.05 |x| \, \mathrm{sgn}(\dot x_1), $$
$$\ddot x_2 = (x_1 - x_2) + (x_3 - x_2) - \dot x_2 - 0.3 \, \mathrm{sgn}(\dot x_2), $$
$$\ddot x_3 = (x_2 - x_3) - \dot x_3 - 0.3 \, \mathrm{sgn} (\dot x_3) + u, $$
or, in a normal form
$$
\dot x_1 = x_4, \ \dot x_2 = x_5, \ \dot x_3 = x_6
$$
$$\dot x_4 = ( - x_1) + (x_2 - x_1) - x_4 - 1.05 |x| \, \mathrm{sgn}(x_4), $$
$$\dot x_5 = (x_1 - x_2) + (x_3 - x_2) - x_5 - 0.3 \, \mathrm{sgn}(x_5), $$
$$\dot x_6 = (x_2 - x_3) - x_6 - 0.3 \, \mathrm{sgn}(x_6) + u, $$
with the initial conditions
$$x_1(0) = -1, \ x_2(0) = 1, \ x_3(0) = -1, $$ $$x_4(0) = -1, \ x_5(0) = 1, \ x_6(0) = 1.$$

\begin{figure*}[h!] \label{pic3}
\begin{minipage}[h]{0.7\linewidth}
\center{\includegraphics[width=1\linewidth]{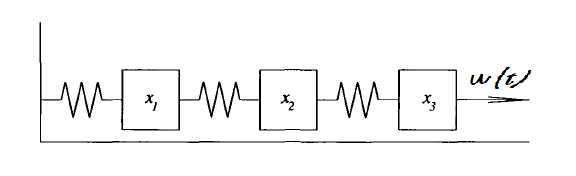} }
\end{minipage}
\caption{Model of Example 7.3}
\label{ris:image2}
\end{figure*}

Due to Filippov definition \cite{fil} the corresponding differential inclusion on the discontinuity surface $x_4 = 0$ is
$$
\dot x_1 = x_4, \ \dot x_2 = x_5, \ \dot x_3 = x_6
$$
$$\dot x_4 \in ( - x_1) + (x_2 - x_1) - x_4 + 1.05 |x| \, [-1, 1] $$
$$\dot x_5 = (x_1 - x_2) + (x_3 - x_2) - x_5 - 0.3 \, \mathrm{sgn}(x_5), $$
$$\dot x_6 = (x_2 - x_3) - x_6 - 0.3 \, \mathrm{sgn}(x_6) + u, $$
with the same boundary conditions.

Suppose that $t^* = -0.2$ and that on the interval $[-0.2, 0]$ no external force is applied, therefore $u(t) = 0$, $t \in [-0.2, 0]$. Such a time interval is enough to bring the system to the vicinity of the surface $x_4(t) = 0$, so we approximately put $x_4(0) = 0$; herewith the other ``initial'' conditions obtained from the movement during the time interval $[-t^*, 0]$ are as follows $x_5(0) \approx 0.026$, $x_6(0) \approx 1.119$, $x_1(0) \approx -1.053$, $x_2(0) \approx 1.099$,  $x_3(0) \approx -0.787$. Take these values as the initial ones for the differential inclusion and apply the paper method in order to find the control $u(t)$ and the correpsonding solution of differential inclusion above on the next time interval $[0, 1]$. The aim of the control $u(t)$ is to bring the third block to the origin at the final time moment $T = 1$, i. e. we impose the restriction $x_3(1) = 0$ on the right endpoint. The ``discontinuous'' surface $x_4 = 0$ of course remains the same.
 
Note that in this example the control $u(t)$ is the function sought as well. The functional structure (put $h_6(x, z, u) = z_6 - ((x_2 - x_3) - x_6 - 0.3 \, \mathrm{sgn}(x_6) + u)$) and the method can be easily modified in a obvious way in order to seek for this variable as well. We also naturally simpify the functions $h_i(x,z)$ and write them down without the $\psi_i$-variable, $i = 1, 2, 3, 5$, since they respond to exact differential equations, not inclusions. Take $(x_{(1)}, z_{(1)}, u_{(1)}) = (0, 0, 0, 0, 0, 0, 0, 0, 0, 0, 0, 0, 0)'$ as the first approximation. In practice we can reduce by $6$ the number of variables (the coordinates $x_i$ and their derivatives $z_i$, $i = \overline{1,3}$), since the first three equations are obviously resolved. At the end of the process the discretization step was equal to $10^{-1}$. Figure 4 illustrates the trajectories obtained. From Figure 4 we see that the relation with the differential inclusion is approximately satisfied. The trajectory practically lies on the surface required as well; therefore the variable $x_1(t)$ remains approximately equal to~$-1.053$ for all $t \in [0, 1]$. The boundary values error doesn't exceed the magnitude $5 \times 10^{-3}$. To obtain such an accuracy $113$ iterations have been required. The functional value on the trajectory obtained is approximately $10^{-3}$. The control obtained is $u_{(113)} = 0.914 t^2 - 4.05 t+ 4.45 t^3$ (we give an approximation via an interpolation polynomial for brevity).

Note that the right-hand side is formally nondifferentiable in phase coordinates function, however it is differentiable if $y_1 \neq 0$ and $y_2 \neq 0$. While algorithm realization $y_1$ and $y_2$ vanishes only at finite number of isolated time moments; hence one can split the time interval in the segments where $y_1$ and $y_2$ remains its signum and Gateaux differentiability of the functional in these variables remains valid.

\begin{figure*}[h!]
\begin{minipage}[h]{0.27\linewidth}
\center{\includegraphics[width=1\linewidth]{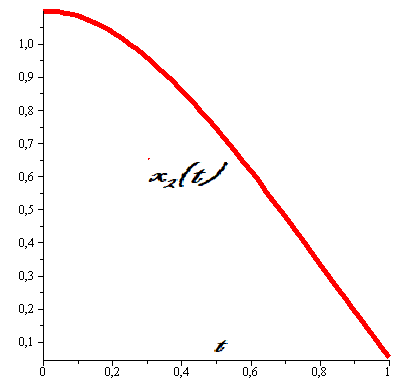} }
\end{minipage}
\hfill
\begin{minipage}[h]{0.27\linewidth}
\center{\includegraphics[width=1\linewidth]{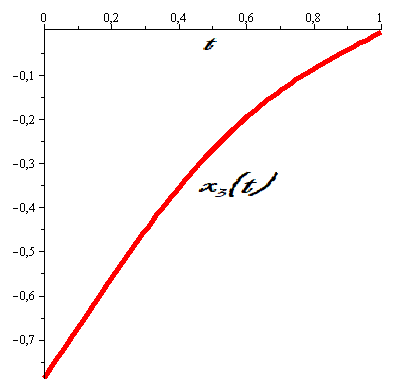} }
\end{minipage}
\hfill
\begin{minipage}[h]{0.27\linewidth}
\center{\includegraphics[width=1\linewidth]{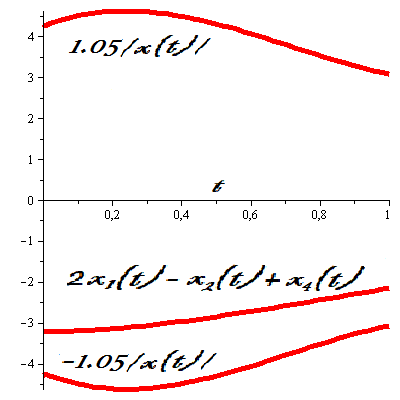} }
\end{minipage}
{\caption{Solution of Example 7.3}}
\end{figure*}
\end{xmpl}

\begin{rmrk}
The detailed discussion on crucial role of an idea to consider $z$ and $\dot x$ as ``independent'' variables in nonsmooth problems of control and variational calculus as well as its general advantages and disadvantages with illustrating examples is given in \cite{Fom_IJC}.
\end{rmrk}


\begin{acknowledgement}{\bf Acknowledgements.}
The author is sincerely greatful to his colleague Maksim Dolgopolik for numerous fruitful discussions and to his friends Tanya Ignatyeva and Dmitriy Kartsev for permanent support. 
\end{acknowledgement}

\end{document}